\RequirePackage{fix-cm}
\documentclass[smallextended]{svjour3}       
\smartqed  
\usepackage{graphicx}
\usepackage{subfigure}
\usepackage{epstopdf}
\usepackage{epsfig}
\DeclareGraphicsExtensions{.pdf, .jpeg, .png}
\usepackage[colorlinks,linkcolor=red]{hyperref}   
\usepackage{bigstrut,multirow,rotating}      
\usepackage{booktabs}                             
\usepackage{enumerate}
\usepackage{amsmath}
\usepackage{amssymb}
\DeclareMathOperator*{\argmin}{arg\,min} 
\usepackage{bbm}
\usepackage{appendix}
\usepackage[numbers,sort&compress]{natbib}
\usepackage{url}
\usepackage[hyphenbreaks]{breakurl}

\begin{document}

\title{Nonlocal Kramers-Moyal formulas and data-driven discovery of stochastic dynamical systems with multiplicative L\'evy noise
}

\titlerunning{Nonlocal Kramers-Moyal formulas}        

\author{Yang Li        \and
        Jinqiao Duan$^*$ 
}


\institute{Li Y. \at
               School of Automation, Nanjing University of Science and Technology, Nanjing 210094, China           
           \and
           Duan J.$^*$ \at
             Department of Mathematics and Department of Physics, Great Bay University, Dongguan, Guangdong 523000, China\\
\email{duan@gbu.edu.cn}
}

\date{Received: date / Accepted: date}

\maketitle

\begin{abstract}
Traditional data-driven methods, effective for deterministic systems or stochastic differential equations (SDEs) with Gaussian noise, fail to handle the discontinuous sample paths and heavy-tailed fluctuations characteristic of L\'evy processes, particularly when the noise is state-dependent. To bridge this gap, we establish nonlocal Kramers-Moyal formulas, rigorously generalizing the classical Kramers-Moyal relations to SDEs with multiplicative L\'evy noise. These formulas provide a direct link between short-time transition probability densities (or sample path statistics) and the underlying SDE coefficients: the drift vector, diffusion matrix, L\'evy jump measure kernel, and L\'evy noise intensity functions. Leveraging these theoretical foundations, we develop novel data-driven algorithms capable of simultaneously identifying all governing components from data and establish convergence results and error analysis for the algorithms. We validate the framework through extensive numerical experiments on prototypical systems. This work provides a principled and practical toolbox for discovering interpretable SDE models governing complex systems influenced by discontinuous, heavy-tailed, state-dependent fluctuations, with broad applicability in climate science, neuroscience, epidemiology, finance, and biological physics.
\end{abstract}

\section{Introduction}
\label{intro}
Environmental fluctuations in complex systems often exhibit statistical properties that cannot be captured by Gaussian processes alone. Non-Gaussian L\'evy motions, particularly $\alpha$-stable L\'evy motions, characterized by heavy-tailed distributions and discontinuous sample paths with large jumps, provide a fundamental mathematical framework for modeling such anomalous dynamics. Rooted in the L\'evy-It\^o decomposition theorem, these processes combine continuous Gaussian diffusion with discontinuous jump components, enabling the representation of a broad spectrum of real-world stochastic perturbations \cite{nunno2008malliavin, applebaum2009levy, Duan2015, bertoin1996levy}. The distinctive properties of L\'evy motion, especially its ability to model extreme events and heavy-tailed distributions, make it indispensable for systems where abrupt, high-amplitude fluctuations significantly influence evolutionary pathways. Consequently, SDEs driven by L\'evy noise have emerged as essential tools for modeling systems influenced by non-Gaussian environmental variability.

Stochastic dynamics induced by L\'evy noise arise ubiquitously across scientific disciplines. In the field of climate science, Ditlevsen found that Earth?s climate variability can be described by a stochastic dynamical system driven by $\alpha$-stable L\'evy noise by analyzing authentic Greenland ice-core measurements \cite{Ditlevsen}. In ecology, Kanazawa et al. have both mathematically and experimentally demonstrated that the motion of tracer particles in active media follows a non-Gaussian L\'evy distribution \cite{kanazawa2020loopy}. Barthelemy et al. experimentally realized L\'evy-flight light diffusion in a tunable optical material, turning L\'evy transport--previously confined to theory and simulation--into an observable and controllable laboratory platform that could enable optical functionalities beyond ordinary diffusion \cite{barthelemy2008levy}. Within quantitative finance, L\'evy noise serves as a potent mathematical construct for characterizing extreme market phenomena such as black swan events \cite{rachev2011financial, tankov2003financial}. Within computational neuroscience, L\'evy noise serves as a critical driver for facilitating stochastic resonance phenomena in neural systems \cite{yamakou2022levy, patel2008stochastic}. Researches suggested that the complex dynamics of electroencephalogram (EEG) signals can be accurately characterized by L\'evy noise \cite{liu2021levy, tajmirriahi2025eeg}. Additionally, recent studies have also pioneered diverse control strategies designed for stochastic systems driven by L\'evy noise, significantly enhancing robustness against jump-diffusion dynamics in complex environments \cite{li2025aperiodically, zhu2024event, zhang2025semi, do2019inverse}. The prevalence of such phenomena underscores the critical need for accurate identification and analysis of SDEs with L\'evy noise from empirical data.

Data-driven methods for discovering dynamical systems have advanced significantly in recent decades, particularly for deterministic models or systems driven by Gaussian noise. Techniques such as Sparse Identification of Nonlinear Dynamics (SINDy) \cite{SINDy1, SINDy2, schaeffer2017learning} and Koopman operator approximations \cite{EDMD, SKO3} excel at extracting explicit governing equations from trajectory data. Similarly, Fokker-Planck-based approaches \cite{hu2025score} and Kramers-Moyal formulas \cite{Boninsegna2018}  have been widely adopted to infer drift and diffusion terms in SDEs with Brownian noise. However, these methods fundamentally rely on local statistical properties and cannot accommodate the discontinuous, heavy-tailed nature of L\'evy motions.

For systems with additive non-Gaussian L\'evy noise, recent progress has been made in estimating stability indices ($\alpha$) and noise intensities from sample path, mean exit times or transition densities \cite{YangLi2020a, liyang2021b, zhang2020extracting, XiaoliChen}. These approaches, however, assume state-independent noise structures, ignoring the multiplicative nature of real-world fluctuations inherent in many jump distributions. Consequently, they fail to address the challenges posed by systems where noise amplitude depends on the system state that restrict their applicability to idealized scenarios.

This work bridges this critical gap by introducing a comprehensive framework for discovering SDEs driven by multiplicative L\'evy noise. Our primary contributions are:
\begin{itemize}
	\item \textbf{Theoretical:} We establish nonlocal Kramers-Moyal formulas, rigorously generalizing the classical Kramers-Moyal relations to SDEs with multiplicative L\'evy noise. These formulas directly link transition probability densities (or sample path statistics) to the drift vector $b(x)$, diffusion matrix $a(x)$, L\'evy kernel $W^{\alpha,\beta}(\xi)$, and state-dependent noise intensity $\sigma(x)$.
	
	\item \textbf{Algorithmic:} Leveraging these formulas, we develop novel data-driven algorithms capable of simultaneously identifying all governing components, i.e., $b(x)$, $a(x)$, $\sigma_i(x_i)$, $\alpha_i$, and $\beta_i$, from sample path data. Crucially, our approach explicitly handles multiplicative noise ($\sigma_i(x_i)$). The convergence results and error analysis for these algorithms are also established and numerically verified.
\end{itemize}

The manuscript is structured as follows: Section \ref{NKMsec} introduces the problem and derives the nonlocal Kramers-Moyal formulas for SDEs with both Brownian noise and L\'evy noise. Section \ref{DDAsec} details the algorithms for learning the L\'evy jump measure, drift coefficient, and diffusion matrix from data and the error analysis. Section \ref{exam} validates the framework through numerical experiments on prototypical systems (Maier-Stein system and R\"ossler oscillator networks). Section \ref{disc} discusses broader implications, limitations, and potential applications.

\section{Nonlocal Kramers-Moyal formulas}
\label{NKMsec}
The Introduction notes that environmental fluctuations often exhibit both Gaussian and non-Gaussian statistical characteristics. This observation aligns with the L\'evy-It\^o decomposition theorem \cite{Duan2015, applebaum2009levy}, which establishes that a wide range of random perturbations can be modeled as linear combinations of a (Gaussian) Brownian motion $B_t$ with continuous sample paths and a (non-Gaussian) pure jump L\'evy motion $L_t$. Therefore, we investigate an $n$-dimensional stochastic dynamical system governed by the following stochastic differential equation
\begin{equation} \label{sde}
{\rm d}x\left( t \right) =b\left( x\left( t- \right) \right){\rm d}t+ \Lambda \left( x\left( t- \right) \right){\rm d}{B_{t}}+ \sigma(x(t-)) {\rm d}{L_{t}}.
\end{equation}
Here ${B_{t}}={{\left[ {{B}_{1,t}},\ \cdots ,\ {{B}_{n,t}} \right]}^{T}}$  and ${L_{t}}={{\left[ {{L}_{1,t}},\ \cdots ,\ {{L}_{n,t}} \right]}^{T}}$ denote $n$-dimensional (Gaussian) Brownian motion and $n$-dimensional (non-Gaussian) L\'evy motion, respectively. This paper focuses primarily on the case where the components of the L\'evy motion are mutually independent. (For details on L\'evy motion, see Appendix \ref{App:A}.) The vector function $b\left( x \right)={{\left[ {{b}_{1}}\left( x \right),\ \cdots ,\ {{b}_{n}} \left( x \right) \right]}^{T}} \in C^1$ represents the drift coefficient (or vector field) in ${{\mathbbm{R}}^{n}}$, while the $n\times n$ matrix $a\left( x \right)= \Lambda(x) {{\Lambda }^{T}(x)} \in C^2$ is the diffusion coefficient. We assume that the $n \times n$ matrix function $\sigma(x)$ is diagonal and continuous, with its diagonal element ${\sigma}_{i}(x)$ representing the noise intensity for the $i$-th component of L\'evy motion. Furthermore, we assume ${\sigma}_{i}(x) > 0$ for $i=1,\dots,n$. (Note that since $-L_{i,t}$ is also a L\'evy motion and the condition ${\sigma}_{i}(x) \neq 0$ is necessary, we strengthen this to positivity for simplicity.) The initial condition of the SDE (\ref{sde}) is $x\left( 0 \right)=z$.

For each component ${{L}_{i,t}}$ ($i=1,...,n$) of the L\'evy motion, its jump measure is assumed to be ${{\nu }_{i}}\left( {\rm d}\xi \right)= {{W}_{i}} \left( \xi \right){\rm d}\xi$ for $\xi\in \mathbbm{R}\backslash \left\{ 0 \right\}$. Among various types of L\'evy motion, we utilize the  $\alpha$-stable L\'evy motion \cite{tankov2003financial, janicki1993simulation} as an example in this paper to design our method, for its broad applicability \cite{OEBbook}. As detailed in Appendix \ref{App:A}, the $\alpha$-stable L\'evy motion is characterized by two key parameters: the stability parameter $\alpha$ (or non-Gaussianity index) and skewness parameter $\beta$ (which governs the motion's symmetry). Denoting ${{W}^{\alpha ,\beta }}\left( \xi \right)$ as the kernel function of the scalar $\alpha$-stable L\'evy motion, it takes the following form:
\begin{equation}\label{kernel}
{{W}^{\alpha ,\beta }}\left( \xi \right)=\left\{ \begin{array}{ccc}
   \frac{{{k}_{\alpha }}\left( 1+\beta  \right)}{2{{\left| \xi \right|}^{1+\alpha }}}, & \xi>0,  \\
   \frac{{{k}_{\alpha }}\left( 1-\beta  \right)}{2{{\left| \xi \right|}^{1+\alpha }}}, & \xi<0,  \\
\end{array} \right.
\end{equation}
where
$$
{{k}_{\alpha }}=\left\{ \begin{array}{ccc}
   \frac{\alpha \left( 1-\alpha  \right)}{\Gamma \left( 2-\alpha  \right)\cos \left( {\pi \alpha }/{2}\; \right)}, & \alpha \ne 1,  \\
   \frac{2}{\pi }, & \alpha =1.  \\
\end{array} \right.
$$

Let $p\left( x,t|z,0 \right)$ denote the transition probability density function of the solution to the stochastic differential equation (\ref{sde}). The corresponding Fokker-Planck equation then takes the following form according to Refs. \cite{Duan2015,dubkov2005generalized}
\begin{equation}\label{FPK}
\begin{aligned}
  \frac{\partial p}{\partial t}= & -\sum\limits_{i=1}^{n} {\frac{\partial } {\partial {{x}_{i}}} \left[ {{b}_{i}} p\left( x,t|z,0 \right) \right]} +\frac{1}{2}\sum\limits_{i,j=1}^{n} {\frac{{{\partial }^{2}}}{\partial {{x}_{i}}\partial {{x}_{j}}} \left[ {{a}_{ij}}p\left( x,t|z,0 \right) \right]} \\
  & +\sum\limits_{i=1}^{n}{\int_{\mathbbm{R}\backslash \left\{ 0 \right\}} {\left[ p\left( x-{{\sigma }_{i}(x)} {y}_{i} {e_{i}}, t|z,0 \right)- p\left( x,t|z,0 \right)+ {{\sigma }_{i}(x)}\chi _{i}^{\alpha }\left( {y}_{i} \right) {y}_{i} \frac{\partial }{\partial {{x}_{i}}} p\left( x,t|z,0 \right) \right] W_{i}^{\alpha ,\beta }\left( {y}_{i} \right){\rm d}{y}_{i}}},
\end{aligned}
\end{equation}
where
$$
\chi _{i}^{\alpha }\left( {y}_{i} \right)=\left\{ \begin{matrix}
   0, & 0<\alpha <1,  \\
   {{\chi }_{\left| \sigma _{i}{y}_{i} \right|<1}}, & \alpha =1,  \\
   1, & 1<\alpha <2.  \\
\end{matrix} \right.
$$
The three terms on the right-hand side of the Fokker-Planck equation correspond to the drift, diffusion, and L\'evy jump components, respectively. Compared to SDEs with pure Gaussian noise, this equation is an integro-partial differential equation featuring a nonlocal term, which can also be expressed as a fractional derivative \cite{Duan2015,dubkov2005generalized}. The symbol ${e_{i}}, i=1,...,n$ in the integral is the standard basis vector in $\mathbbm{R}^n$, and ${{\chi }_{\left| \sigma _{i}{y}_{i} \right|<1}}$ is the indicator function. The initial condition for the Fokker-Planck equation is $p\left( x,0|z,0 \right)=\delta \left( x-z \right)$. Additionally, we define the $n-1$-dimensional vector $x^i$ as the vector $x$ with the $i$-th component $x_i$ removed, and $p_i(x_i,t|z,0) =\int_{\mathbbm{R}^{n-1}} p(x,t|z,0){\rm d}x^i$ as the marginal probability distribution along $i$-th direction. The set $\Gamma $ in the following theorem is defined as an $n$-dimensional cube ${{\left[ -\varepsilon ,\ \varepsilon  \right]}^{n}}$.

This paper aims to infer stochastic dynamical systems driven by multiplicative (Gaussian) Brownian noise and multiplicative (non-Gaussian) L\'evy noise from sample path data. To this end, we first establish the following theorem, which expresses the drift coefficient, diffusion coefficient, L\'evy jump measure and noise intensity in terms of the probability density function $p\left( x,t|z,0 \right)$ (the solution to the Fokker-Planck equation (\ref{FPK})).

\newtheorem{thm}{\bf Theorem}
\begin{thm} (Relation between stochastic governing law and Fokker-Planck equation)\\
\label{thm1}
For every $\varepsilon >0$, the probability density function $p\left( x,t|z,0 \right)$ and the jump measure, drift and diffusion have the following relations:
\begin{enumerate}[1)]
\item For every ${{x}_{i}}$ and $z$ satisfying $\left| {{x}_{i}}-{{z}_{i}} \right|>\varepsilon $ and $i=1,\ 2,\ \ldots ,\ n$,
$$
\underset{t\to 0}{\mathop{\lim }}\,{{t}^{-1}}{{p}_{i}} \left( {{x}_{i}},t|z,0 \right)= \sigma _{i}^{-1}(z) W_{i}^{\alpha ,\beta } \left( \sigma _{i}^{-1}(z) \left( {{x}_{i}}-{{z}_{i}} \right) \right).
$$

\item For $i=1,\ 2,\ \ldots ,\ n$,
$$
\underset{t\to 0}{\mathop{\lim }}\,{{t}^{-1}}\int_{x-z\in \Gamma }{\left( {{x}_{i}}-{{z}_{i}} \right) p\left( x,t|z,0 \right) {\rm d}x} ={{b}_{i}} \left( z \right)+R_{i}^{\alpha ,\beta, \varepsilon}\left( z \right),
$$
where $R_{i}^{\alpha ,\beta, \varepsilon}\left( z \right)=\left\{ \begin{array}{lc}
   \sigma _{i}^{-1}(z)\int_{-\varepsilon }^{\varepsilon }{{y}_{i}W_{i}^{\alpha ,\beta } \left( \sigma _{i}^{-1}(z){y}_{i} \right) {\rm d}{y}_{i}}, & \alpha <1,  \\
   \sigma _{i}^{-1}(z)\left[ \int_{-\varepsilon }^{-1}{{y}_{i}W_{i}^{\alpha ,\beta }\left( \sigma _{i}^{-1}(z){y}_{i}\right) {\rm d}{y}_{i}}+ \int_{1}^{\varepsilon }{{y}_{i}W_{i}^{\alpha ,\beta }\left( \sigma _{i}^{-1}(z){y}_{i} \right){\rm d}{y}_{i}} \right], & \alpha =1,  \\
   -\sigma _{i}^{-1}(z)\left[ \int_{-\infty }^{-\varepsilon } {{y}_{i} W_{i}^{\alpha ,\beta } \left( \sigma _{i}^{-1}(z){y}_{i} \right) {\rm d} {y}_{i}}+ \int_{\varepsilon }^{\infty }{{y}_{i}W_{i}^{\alpha ,\beta }\left( \sigma _{i}^{-1}(z){y}_{i} \right){\rm d}{y}_{i}} \right], & \alpha >1.  \\
\end{array} \right.$

\item For $i,j=1,\ 2,\ \ldots ,\ n$,
$$
\underset{t\to 0}{\mathop{\lim }}\,{{t}^{-1}}\int_{x-z\in \Gamma }{\left( {{x}_{i}}-{{z}_{i}} \right)\left( {{x}_{j}}-{{z}_{j}} \right)p\left( x,t|z,0 \right){\rm d}x}={{a}_{ij}}\left( z \right)+ S_{ij}^{\alpha ,\beta, \varepsilon} \left( z  \right),
$$
where $S_{ii}^{\alpha ,\beta, \varepsilon} \left( z  \right)=\sigma _{i}^{-1}(z) \int_{-\varepsilon }^{\varepsilon }{{{y}_{i}^{2}}W_{i}^{\alpha ,\beta }\left( \sigma _{i}^{-1}(z){y}_{i} \right){\rm d}{y}_{i}}$ and $S_{ij}^{\alpha ,\beta, \varepsilon} \left( z  \right)=0$ for $i \neq j$.
\end{enumerate}
\end{thm}

Note that the core idea of Theorem \ref{thm1} originates from the following observation: sample paths of SDEs with Gaussian noise exhibit continuity, whereas those of L\'evy noise exhibit discontinuities. Consequently, large jumps occurring within short time intervals must originate from L\'evy noise. Based on this principle, we partition the state space to isolate contributions from L\'evy noise, drift, and diffusion. However, designing numerical algorithms directly using probability density functions is computationally challenging. Therefore, we restate Theorem \ref{thm1} as Corollary \ref{cor2} below, which characterizes the system's drift, diffusion, and L\'evy jump measure in terms of sample paths (specifically, solutions to the SDE (\ref{sde})).

\newtheorem{cor}[thm]{\bf Corollary} 
\begin{cor}(Nonlocal Kramers-Moyal formulas)\\
\label{cor2}
For every $\varepsilon >0$, the sample path solution $x\left( t \right)$ of the stochastic differential equation (\ref{sde}) and the jump measure, drift and diffusion have the following relations:
\begin{enumerate}[1)]
\item For every ${{c}_{1}}$ and ${{c}_{2}}$ satisfying ${{c}_{1}}<{{c}_{2}}<0$ or $0<{{c}_{1}}<{{c}_{2}}$, and $i=1,\ 2,\ \ldots ,\ n$,
$$
\underset{t\to 0}{\mathop{\lim }}\,{{t}^{-1}}\mathbbm{P}\left\{ \left. {{x}_{i}} \left( t \right)-{{z}_{i}}\in \left[ {{c}_{1}}, {{c}_{2}} \right) \right|x\left( 0 \right)=z \right\}=\sigma _{i}^{-1}(z) \int_{{{c}_{1}}}^{{{c}_{2}}} {W_{i}^{\alpha ,\beta }\left( \sigma _{i}^{-1}(z) {y}_{i} \right){\rm d}{y}_{i}}.
$$
\item For $i=1,\ 2,\ \ldots ,\ n$,
$$
\begin{aligned}
  & \underset{t\to 0}{\mathop{\lim }}\,{{t}^{-1}}\mathbbm{P}\left\{ \left. x\left( t \right)-z\in \Gamma  \right|x\left( 0 \right)=z \right\}\cdot \mathbbm{E}\left[ \left. \left( {{x}_{i}}\left( t \right)-{{z}_{i}} \right) \right|x\left( 0 \right)=z;\ x\left( t \right)-z\in \Gamma  \right] \\
 & ={{b}_{i}}\left( z \right)+R_{i}^{\alpha ,\beta, \varepsilon}\left( z \right). \\
\end{aligned}
$$
\item For $i,j=1,\ 2,\ \ldots ,\ n$,
$$
\begin{aligned}
  & \underset{t\to 0}{\mathop{\lim }}\,{{t}^{-1}}\mathbbm{P}\left\{ \left. x\left( t \right)-z\in \Gamma  \right|x\left( 0 \right)=z \right\}\cdot \mathbbm{E}\left[ \left. \left( {{x}_{i}}\left( t \right)-{{z}_{i}} \right)\left( {{x}_{j}}\left( t \right)-{{z}_{j}} \right) \right|x\left( 0 \right)=z;\ x\left( t \right)-z\in \Gamma  \right] \\
 & ={{a}_{ij}}\left( z \right)+S_{ij}^{\alpha ,\beta, \varepsilon} \left( z  \right). \\
\end{aligned}
$$
\end{enumerate}
\end{cor}

Note that the three formulas in this corollary establish a direct connection between sample paths and the coefficient functions of SDE (\ref{sde}). These can be viewed as generalizations of the Kramers-Moyal formulas to systems with non-Gaussian noise. Consequently, we term them nonlocal Kramers-Moyal formulas, contrasting with the standard (local) Kramers-Moyal formulas for SDEs with Gaussian noise as presented in \cite[Ch. 3]{KM}.

Proofs of Theorem \ref{thm1} and Corollary \ref{cor2} appear in Appendices \ref{App:B} and \ref{App:C}, respectively.

Building on this formulation, we derive nonlocal Kramers-Moyal formulas for SDEs of the form (\ref{sde}) with asymmetric L\'evy noise whose components are mutually independent. However, rotationally symmetric L\'evy noise which is another common type lies beyond this framework due to inherent cross-component correlations. While this work focuses primarily on the case described by Eq. (\ref{sde}), we present corresponding theoretical results for the rotationally symmetric scenario to ensure completeness in Appendix \ref{App:D}.

\section{Methods}
\label{DDAsec}
In this section, we propose data-driven modeling algorithms for learning stochastic dynamical systems with multiplicative (Gaussian) Brownian noise and multiplicative (non-Gaussian) L\'evy noise from sample path data. We focus primarily on designing algorithms for SDEs of the form (\ref{sde}), and will discuss the rorationally symmetric case in Appendix \ref{App:D}.

\begin{figure}
	\centering
	\includegraphics[width=12cm]{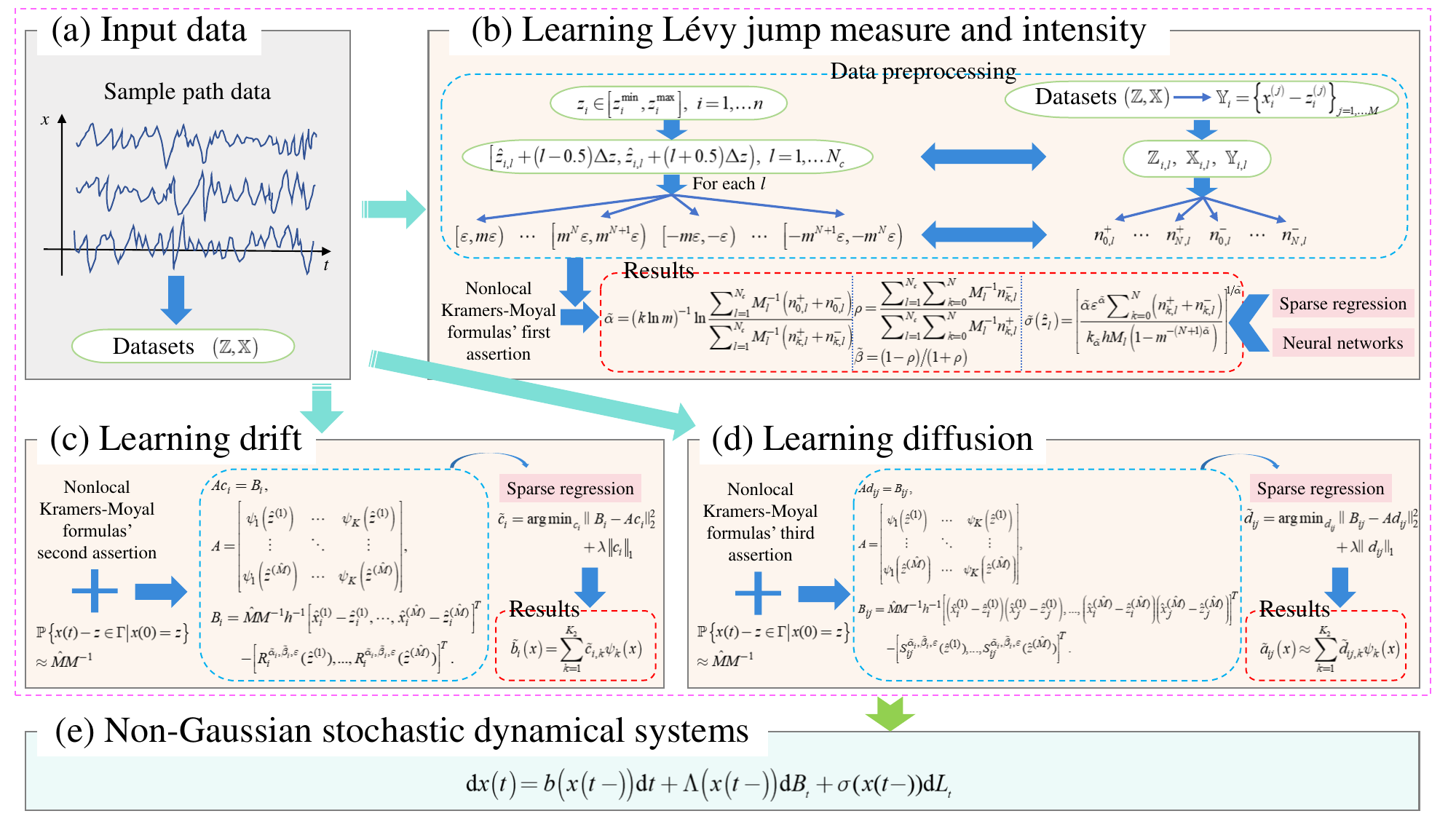}
	\caption{Schematic diagrams for illustrating the proposed data-driven framework. (a) Input data, including sample path data and the constructed datasets $(\mathbbm{Z},\mathbbm{X})$. (b) The algorithm for learning L\'evy jump measure and intensity functions, including data preprocessing and L\'evy noise parameter identification such as non-Gaussian index $\alpha$, skewness parameter $\beta$, and noise intensity function $\sigma(x)$. (c) The algorithm for learning the drift coefficient $b(x)$. (d) The algorithm for learning the diffusion coefficient $a(x)$. (e) The underestimated non-Gaussian stochastic dynamical systems.}
	\label{fig4}
\end{figure}

Consider paired datasets $\mathbbm{Z}$ and $\mathbbm{X}$
\begin{equation}\label{eq3}
\begin{aligned}
  & \mathbbm{Z}=\left[ {z^{(1)}},\ {z^{(2)}},\ \cdots ,\ {z^{(M)}} \right], \\
  & \mathbbm{X}=\left[ {x^{(1)}},\ {x^{(2)}},\ \cdots ,\ {x^{(M)}} \right], \\
\end{aligned}
\end{equation}
where each ${x^{(j)}} \in \mathbbm{R}^n$ is the image of ${z^{(j)}} \in \mathbbm{R}^n$ after a time step $h$. We then devise algorithms to discover the L\'evy motion's kernel function and noise intensity, the drift coefficient, and the diffusion matrix from these datasets in the following three subsections. To enhance the understanding of the algorithmic structure, we present the complete architecture of the algorithms in Fig. \ref{fig4}.

\subsection{Algorithm for learning the L\'evy motion}
\label{learnlevysubsec}
To facilitate algorithm design and prevent the curse of dimensionality, we assume the L\'evy noise intensity satisfies $\sigma_i(x)=\sigma_i(x_i)$, $i=1,\ 2,\ \ldots ,\ n$. That is, every $\sigma_i(x)$ is a scalar function depending only on its corresponding coordinate $x_i$. We design algorithms for this $\sigma_i(x_i)$ case and will discuss the general $\sigma_i(x)$ scenario later.

According to the kernel function expression (\ref{kernel}), determining the L\'evy motion requires learning only the parameters ${{\alpha }_{i}}$, ${{\beta }_{i}}$ and noise intensity functions ${{\sigma }_{i}}(x_i)$ for $i=1,\ 2,\ \ldots ,\ n$. Corollary \ref{cor2}'s first assertion shows that the equation depends solely $z_i$ and $x_i-z_i$, independent of the other $n-1$ dimensions. This allows separate consideration of each component. We define $\mathbbm{Z}_{i}$ and $\mathbbm{X}_{i}$ as the collections of the $i$-th components from $\mathbbm{Z}$ and $\mathbbm{X}$, respectively, and construct new datasets ${\mathbbm{Y}_{i}}= \left[ y_{i}^{(1)},\ y_{i}^{(2)},\ \cdots ,\ y_{i}^{(M)} \right]$ where ${{y}_{i}^{(j)}}={{x}_{i}^{(j)}}-{{z}_{i}^{(j)}}$ for $i=1,\ 2,\ \ldots ,\ n$, $j=1,\ 2,\ \ldots ,\ M$.

Since Corollary \ref{cor2}'s first assertion depends on $z_i$, we partition its phase space $[z_i^{\text{min}},z_i^{\text{max}}]$ into $N_c$ intervals: $\left[ z_i^{\text{min}},\ z_i^{\text{min}}+ \Delta z_i \right)$, $\left[ z_i^{\text{min}}+ \Delta z_i,\ z_i^{\text{min}}+ 2\Delta z_i \right)$, $\ldots$, $\left[ z_i^{\text{max}}- \Delta z_i,\ z_i^{\text{max}} \right)$, where $\Delta z_i= (z_i^{\text{max}}-z_i^{\text{min}})/N_c$. Each interval's midpoint $\hat{z}_{i,l}= z_i^{\text{min}}+ (l-0.5) \Delta z_i$, $l=1,\ 2,\ \ldots ,\ N_c$ enables approximation of ${{\sigma}_{i}}(z_i)$ as constant ${{\sigma }_{i}}(\hat{z}_{i,l})$. Within each interval $\left[ z_i^{\text{min}}+ (l-1) \Delta z_i,\ z_i^{\text{min}}+ l\Delta z_i \right)$, we estimate ${{\alpha }_{i}}$, ${{\beta }_{i}}$ and ${{\sigma }_{i}}(\hat{z}_{i,l})$ using Corollary \ref{cor2}'s first assertion. Let $\mathbbm{Z}_{i,l}$ denote $\mathbbm{Z}_i$'s subset in this interval, with $\mathbbm{X}_{i,l}$ and $\mathbbm{Y}_{i,l}$ being the corresponding $M_l$-element datasets. The conditional probability is then estimated as the ratio of points in $\mathbbm{Y}_{i,l}$ falling within $\left[ {{c}_{1}},\ {{c}_{2}} \right)$ to $M_l$.

Specifically, we choose $2N+2$ intervals $\left[ \varepsilon ,\ m\varepsilon  \right)$, $\left[ m\varepsilon ,\ {{m}^{2}}\varepsilon  \right)$, $\ldots$, $\left[ {{m}^{N}}\varepsilon ,\ {{m}^{N+1}}\varepsilon  \right)$ and $\left[ -m\varepsilon ,\ -\varepsilon  \right)$, $\left[ -{{m}^{2}}\varepsilon ,\ -m\varepsilon  \right)$, $\ldots$, $\left[ -{{m}^{N+1}}\varepsilon ,\ -{{m}^{N}}\varepsilon  \right)$ with the positive integer $N$, the positive real number $\varepsilon $, and the real number $m>1$. Suppose that there are $n_{0,l}^{+},\ n_{1,l}^{+},\ \ldots ,\ n_{N,l}^{+}$ and $n_{0,l}^{-},\ n_{1,l}^{-},\ \ldots ,\ n_{N,l}^{-}$ points from the data set $\mathbbm{Y}_{i,l}$ falling into these intervals. For convenience, we omit the dependence of $n_{k,l}^{+(-)}$ about the index $i$ and also denote ${{\alpha }_{i}}$, ${{\beta }_{i}}$, ${{\sigma }_{i}}(\hat{z}_{i,l})$ as $\alpha$, $\beta$, $\sigma_l$ in the following since we can deal with every component of the L\'evy motion independently.

From Corollary \ref{cor2}'s first assertion:
$$
\begin{aligned}
& {{h}^{-1}}\mathbbm{P}\left\{ \left. {{x}_{i}}\left( h \right)-{{z}_{i}}\in \left[ {{m}^{k}}\varepsilon ,\ {{m}^{k+1}}\varepsilon  \right) \right|x\left( 0 \right)=z \right\}\approx {{h}^{-1}}{{M}_l^{-1}}n_{k,l}^{+}, \\
& k=0,\ 1,\ \ldots ,\ N, \ l=1,\ 2,\ \ldots,\ N_c.
\end{aligned}
$$
Meanwhile, its right-hand side integration can be expressed as
$$
\begin{aligned}
  & {{\sigma }_l^{-1}}\int_{{{m}^{k}}\varepsilon }^{{{m}^{k+1}}\varepsilon }{W_{i}^{\alpha ,\beta }\left( {{\sigma }_l^{-1}}{y}_{i} \right){\rm d}{y}_{i}} \\
  & ={{{\sigma }_l^{\alpha }}{{k}_{\alpha }}}\frac{\left( 1+\beta  \right)}{2}\int_{{{m}^{k}}\varepsilon }^{{{m}^{k+1}}\varepsilon }{{{\left| {y}_{i} \right|}^{-\left( 1+\alpha  \right)}}{\rm d}{y}_{i}} \\
  & ={{{\sigma }_l^{\alpha }}{{k}_{\alpha }}\left( 1+\beta  \right){{\alpha }^{-1}}{{\varepsilon }^{-\alpha }}{{m}^{-k\alpha }}\left( 1-{{m}^{-\alpha }} \right)}/{2}, \\
  & k=0,\ 1,\ \ldots ,\ N, \ l=1,\ 2,\ \ldots,\ N_c.
\end{aligned}
$$
Combining these two equations yields $(N+1)N_c$ equalities
\begin{equation}\label{eq4}
	\begin{aligned}
    & {{{\sigma }_l^{\alpha }}{{k}_{\alpha }}\left( 1+\beta  \right){{\alpha }^{-1}}{{\varepsilon }^{-\alpha }}{{m}^{-k\alpha }}\left( 1-{{m}^{-\alpha }} \right)}/{2}={{h}^{-1}}{{M}_l^{-1}}n_{k,l}^{+}, \\
    & k=0,\ 1,\ \ldots ,\ N, \ l=1,\ 2,\ \ldots,\ N_c.
    \end{aligned}
\end{equation}
On the other $(N+1)N_c$ negative intervals, we similarly obtain
\begin{equation}\label{eq5}
\begin{aligned}
	& {{{\sigma }_l^{\alpha }}{{k}_{\alpha }}\left( 1-\beta  \right){{\alpha }^{-1}}{{\varepsilon }^{-\alpha }}{{m}^{-k\alpha }}\left( 1-{{m}^{-\alpha }} \right)}/{2}={{h}^{-1}}{{M}_l^{-1}}n_{k,l}^{-}, \\
	& k=0,\ 1,\ \ldots ,\ N, \ l=1,\ 2,\ \ldots,\ N_c.
\end{aligned}
\end{equation}
Summing Eqs. (\ref{eq4}) and (\ref{eq5}) produces
\begin{equation}\label{eq5a}
\begin{aligned}
	& {{{\sigma }_l^{\alpha }}{{k}_{\alpha }}{{\alpha }^{-1}}{{\varepsilon }^{-\alpha }}{{m}^{-k\alpha }}\left( 1-{{m}^{-\alpha }} \right)}={{h}^{-1}}{{M}_l^{-1}}(n_{k,l}^{+}+n_{k,l}^{-}), \\
	& k=0,\ 1,\ \ldots ,\ N, \ l=1,\ 2,\ \ldots,\ N_c.
\end{aligned}
\end{equation}
Summing Eq. (\ref{eq5a}) over index $l$ gives
\begin{equation}\label{eq5b}
	\begin{aligned}
		& {{{k}_{\alpha }}{{\alpha }^{-1}}{{\varepsilon }^{-\alpha }} {{m}^{-k\alpha }}\left( 1-{{m}^{-\alpha }} \right)} \sum_{l=1}^{N_c} {{\sigma }_l^{\alpha }} ={{h}^{-1}}\sum_{l=1}^{N_c} {{M}_l^{-1}} (n_{k,l}^{+}+n_{k,l}^{-}), \\
		& k=0,\ 1,\ \ldots ,\ N.
	\end{aligned}
\end{equation}
Taking ratios of the first equation ($k=0$) to the remaining $N$ equations yields
\begin{equation}\label{eq6}
\alpha ={{\left( k\ln m \right)}^{-1}}\ln \frac{\sum_{l=1}^{N_c} {{M}_l^{-1}} (n_{0,l}^{+}+n_{0,l}^{-})}{\sum_{l=1}^{N_c} {{M}_l^{-1}} (n_{k,l}^{+}+n_{k,l}^{-})},\quad k=1,\ 2,\ \ldots ,\ N.
\end{equation}
For $N=1$, this provides the optimal estimate $\tilde{\alpha}$. When $N \geq 2$, we take $\tilde{\alpha}$ as the mean value of Eqs. (\ref{eq6}).

Let
\begin{equation}\label{eq7a}
\rho =\frac{\sum_{l=1}^{N_c} \sum_{k=0}^{N} {{M}_l^{-1}} {n_{k,l}^{-}}} {\sum_{l=1}^{N_c} \sum_{k=0}^{N} {{M}_l^{-1}} {n_{k,l}^{+}}}.
\end{equation}
The parameter $\beta$ is then determined by $\rho$, which equals the ratio of the double summation (over $k$ and $l$) of Eq. (\ref{eq4}) to that of Eq. (\ref{eq5})
\begin{equation}\label{eq7}
\beta =\frac{1-\rho }{1+\rho }.
\end{equation}

Summing Eq. (\ref{eq5a}) over index $k$ yields
\begin{equation}\label{eq8a}
	\begin{aligned}
		& {{{\sigma }_l^{\alpha }}{{k}_{\alpha }}{{\alpha }^{-1}}{{\varepsilon }^{-\alpha }}\left( 1-{{m}^{-(N+1)\alpha }} \right)}= {{h}^{-1}} {{M}_l^{-1}} \sum_{k=0}^{N}(n_{k,l}^{+}+n_{k,l}^{-}), \\
		& l=1,\ 2,\ \ldots,\ N_c.
	\end{aligned}
\end{equation}
Thus, the noise intensity $\sigma_l$ can be approximated as
\begin{equation}\label{eq8}
\sigma_l ={{\left[ \frac{\tilde{\alpha }{{\varepsilon }^{{\tilde{\alpha }}}} \sum_{k=0}^{N}(n_{k,l}^{+}+n_{k,l}^{-})} {{{k}_{{ \tilde{\alpha }}}} hM_l \left( 1-{{m}^{-(N+1)\tilde{\alpha } }} \right)} \right]} ^{{1}/{{\tilde{\alpha }}}\;}},\quad l=1,\ 2,\ \ldots,\ N_c.
\end{equation}

After applying Eqs. (\ref{eq6}), (\ref{eq7}), and (\ref{eq8}), we obtain the approximation parameters ${\tilde{\alpha }_{i}}$, ${\tilde{\beta }_{i}}$ and noise intensity function values ${\tilde{\sigma }_{i}}(\hat{z}_{i,l})$ at $N_c$ points. Processing all $n$ dimensions yields the complete parameter set and functions needed to determine both the L\'evy jump measure and noise intensity.

Define $\mathbbm{Q}_i= \{ {\tilde{\sigma }_{i}}(\hat{z}_{i,l}): l=1, 2, \ldots,\ N_c \}$, $i=1, 2, ..., n$. Regression methods then yield approximation functions $\sigma_{i,\theta}(x_i)$ using training set $\mathbbm{Q}_i$, where $\theta$ represents the model's training parameters. This regression problem minimizes the loss function
\begin{equation}\label{loss}
\mathcal{L}_i= \frac{1}{N_c} \sum_{l=1}^{N_c} | \sigma_{i,\theta}(\hat{z}_{i,l})- {\tilde{\sigma }_{i}}(\hat{z}_{i,l}) |^2.
\end{equation}
For example, fully connected neural networks can approximate the noise intensity functions. Each dimension's network has input $x_i$ and output $\sigma_{i,\theta}$, with $\theta$ denoting all weights and biases.

Sparse regression offers another approach for estimating noise intensity functions. Consider a dictionary of candidate functions $\Phi(\xi) = \{ \phi_1(\xi), \phi_2(\xi), ..., \phi_{K_1}(\xi) \}$ where $\xi \in \mathbbm{R}$. The approximation then takes the form $\sigma_{i,\theta}(x_i) = \theta_1 \phi_1(x_i)+ \theta_2 \phi_2(x_i)+ \cdots +\theta_{K_1} \phi_{K_1}(x_i)$. Let
\begin{equation}\label{srAB}
	A_i= \begin{bmatrix}
		\phi_1(\hat{z}_{i,1}) & \phi_2(\hat{z}_{i,1}) & \cdots & \phi_{K_1}(\hat{z}_{i,1}) \\ \phi_1(\hat{z}_{i,2}) & \phi_2(\hat{z}_{i,2}) & \cdots & \phi_{K_1}(\hat{z}_{i,2}) \\ \vdots & \vdots & \ddots & \vdots \\ \phi_1(\hat{z}_{i,N_c}) & \phi_2(\hat{z}_{i,N_c}) & \cdots & \phi_{K_1}(\hat{z}_{i,N_c})
	\end{bmatrix},
    \ B_i= \begin{bmatrix}
	    {\tilde{\sigma }_{i}}(\hat{z}_{i,1}) \\ {\tilde{\sigma }_{i}}(\hat{z}_{i,2}) \\ \vdots \\ {\tilde{\sigma }_{i}}(\hat{z}_{i,N_c})
    \end{bmatrix}.
\end{equation}
Consequently, minimizing the loss function (\ref{loss}) reduces to the regression problem
\begin{equation}\label{sigmaloss1}
	\tilde{\theta}= \argmin_{\theta} \| B_i- A_i \theta \|_2^2,
\end{equation}
where $\left\| \cdot \right\|_2$ denotes the 2-norm. To promote sparsity, we augment Eq. (\ref{sigmaloss1}) with a regularization term
\begin{equation}\label{sigmaloss2}
	\tilde{\theta}= \argmin_{\theta} \| B_i- A_i \theta \|_2^2 +\lambda \| \theta \|_1,
\end{equation}
with $\left\| \cdot \right\|_1$ denoting the 1-norm. This sparse regression problem (\ref{sigmaloss2}) can be efficiently solved via either the SINDy method or stepwise sparse regressor algorithm, detailed in Appendices \ref{App:E} and \ref{App:G}, respectively.

Several remarks apply to this algorithm. First, symmetric L\'evy motion components yield $n_{k}^{+} = n_{k}^{-}$, implying $\rho =1$ and $\beta =0$ in Eqs. (\ref{eq7a}) and (\ref{eq7}) which are consistent with symmetric behavior. Second, different $N_c$ choices (e.g., 10, 11) generate corresponding numbers of $\sigma_i$ function values, enabling denser sampling through multiple $N_c$ selections. Third, crucially, this method avoids the curse of dimensionality by processing each component independently: dimensionality increases do not necessarily demand more data. This advantage stems from assuming $\sigma_i(x_i)$ depends solely on $x_i$. Should $\sigma_i$ depend on all state variables, separate dimensional processing becomes impossible, requiring full phase space bin discretization, an approach demanding exponentially more data with increasing dimension.

\subsection{Algorithm for learning the drift coefficient}
\label{learndriftsubsec}
Based on the second formula of Corollary \ref{cor2}, we design a data-driven algorithm to learn the drift coefficient $b(x)$ from datasets $\mathbbm{Z}$ and $\mathbbm{X}$. Assume another dictionary of candidate functions $\Psi(x) = \{ \psi_1(x), \psi_2(x), ..., \psi_{K_2}(x) \}$ where $x \in \mathbbm{R}^n$. Each component of the drift term can be approximated as ${{b}_{i}}\left( x \right)\approx \sum_{k=1}^{K_2}{{{c}_{i,k}} {{\psi }_{k}} \left( x \right)}$, $i=1, 2, ..., n$, where $c_i= [c_{i,1}, c_{i,2}, ..., c_{i,K_2}]^T$ denotes the coefficient vector for the $i$-th component. Therefore, the sparse regression method can also be used to learn the drift coefficient.

Note that substituting the kernel function (\ref{kernel}) of the L\'evy jump measure into the expression for $R_{i}^{\alpha ,\beta, \varepsilon}(z)$ from Corollary \ref{cor2}'s second assertion yields
\begin{equation}\label{Rz1}
	R^{\alpha ,\beta, \varepsilon}(z)= \left\{
	\begin{array}{ll}
		\sigma^{\alpha}(z) k_{\alpha} \beta \frac{\varepsilon^{1-\alpha}} {1-\alpha}, \ \ & \alpha \neq 1, \\
		\frac{2}{\pi} \sigma(z) \beta \ln\varepsilon, & \alpha =1,
	\end{array} \right.
\end{equation}
where we omit the dependence of $\alpha$, $\beta$, $\sigma$ on dimension $i$ for simplicity. Therefore, $R^{\alpha ,\beta, \varepsilon}(z)$ can be computed numerically using the estimated parameters $\tilde{\alpha}_i$, $\tilde{\beta}_i$ and the learned noise intensity function $\sigma_{i,\theta}(x_i)$ from the preceding subsection.

According to Corollary \ref{cor2}'s second assertion, the drift term estimation depends only on data satisfying $x(h)-z \in \Gamma$. After removing data outside the cube $\Gamma$ from datasets $\mathbbm{Z}$ and $\mathbbm{X}$, we obtain new datasets with $\hat{M}$ elements each
\begin{equation}\label{datasets2}
	\begin{aligned}
		& \mathbbm{\hat{Z}}=\left[ {\hat{z}^{(1)}},\ {\hat{z}^{(2)}},\ \cdots ,\ {\hat{z}^{(\hat{M})}} \right], \\
		& \mathbbm{\hat{X}}=\left[ {\hat{x}^{(1)}},\ {\hat{x}^{(2)}},\ \cdots ,\ {\hat{x}^{(\hat{M})}} \right]. \\
	\end{aligned}
\end{equation}
Note that the conditional probability $\mathbbm{P} \left\{ \left. x(t)-z\in \Gamma  \right| x(0)=z \right\}$ in Corollary \ref{cor2}'s second formula depends on the initial position $z$. However, under the condition $\sqrt{h} \ll \varepsilon$ where large jumps originate solely from the L\'evy motion, we have $\mathbbm{P} \left\{ \left. x(h)-z \notin \Gamma  \right| x(0)=z \right\} \approx \mathbbm{P} \left\{ x(h)-z \notin \Gamma \right\}$. This implies $\mathbbm{P} \left\{ \left. x(h)-z \in \Gamma  \right| x(0)=z \right\} \approx \mathbbm{P} \left\{ x(h)-z \in \Gamma \right\}$. Therefore, this conditional probability can be numerically approximated as $\mathbbm{P} \left\{ \left. x(t)-z \in \Gamma  \right| x(0)=z \right\} \approx \frac{\hat{M}}{M}$.

Assume the system is ergodic, allowing the conditional expectations in the nonlocal Kramers-Moyal formulas to be replaced by a single trajectory. Consequently, we derive the following system of linear equations
\begin{equation}\label{eqdrift}
\begin{aligned}
  & A{c_{i}}={{B}_{i}}, \\
 & A=\left[ \begin{array}{ccc}
   {{\psi }_{1}}\left( \hat{z}^{(1)} \right) & \cdots  & {{\psi }_{K}}\left( \hat{z}^{(1)} \right)  \\
   \vdots  & \ddots  & \vdots   \\
   {{\psi }_{1}}\left( \hat{z}^{(\hat{M})} \right) & \cdots  & {{\psi }_{K}}\left( \hat{z}^{(\hat{M})} \right)  \\
\end{array} \right], \\
 & {{B}_{i}}= \hat{M}{{M}^{-1}} {{h}^{-1}}{{\left[ \hat{x}_i^{(1)}- \hat{z}_i^{(1)},\ \cdots ,\ \hat{x}_i^{(\hat{M})}-\hat{z}_i^{(\hat{M})} \right]}^{T}} \\
 & \ \ \ \ \ \ \ - \left[ R_{i}^{\tilde{\alpha}_i, \tilde{\beta}_i, \varepsilon} (\hat{z}^{(1)}), ..., R_{i}^{\tilde{\alpha}_i, \tilde{\beta}_i, \varepsilon} (\hat{z}^{(\hat{M})}) \right]^T. \\
\end{aligned}
\end{equation}
This system can be transformed into a least squares regression problem
\begin{equation}\label{driftloss1}
	\tilde{c}_i= \argmin_{c_i} \| B_i- A c_i \|_2^2.
\end{equation}
To promote sparsity, we add a regularization term to Eq. (\ref{driftloss1}), yielding
\begin{equation}\label{driftloss2}
	\tilde{c}_i= \argmin_{c_i} \| B_i- A c_i \|_2^2 +\lambda \| c_i \|_1.
\end{equation}

SINDy, introduced in Appendix \ref{App:E}, is a classical method for solving the sparse regression problem (\ref{driftloss2}). For low-dimensional systems ($n=1,2,3$), a combination of binning, stepwise sparse regression, and cross-validation algorithms provides a more accurate, parameter-free alternative, as detailed in Appendices \ref{App:F} and \ref{App:G}. However, binning is unsuitable for high-dimensional systems. Thus, we employ the combined method for low-dimensional systems and SINDy for high-dimensional cases. After processing all $n$ dimensions (indexed by $i$), we successfully learn the drift vector function from data.

\subsection{Algorithm for learning the diffusion coefficient}
\label{learndiffusionsubsec}
Based on Corollary \ref{cor2}'s third formula, we design a data-driven algorithm to learn the diffusion coefficient $a(x)$ from datasets $\mathbbm{Z}$ and $\mathbbm{X}$, analogous to the drift term learning algorithm. Using the same dictionary of candidate functions $\Psi(x) = \{ \psi_1(x), \psi_2(x), ..., \psi_{K_2}(x) \}$, $x \in \mathbbm{R}^n$, we approximate each component of the diffusion term as ${{a}_{ij}} \left( x \right) \approx \sum_{k=1}^{K_2}{{{d}_{ij,k}} {{\psi }_{k}} \left( x \right)}$, $i,j=1, 2, ..., n$, where $d_{ij}= [d_{ij,1}, d_{ij,2}, ..., d_{ij,K_2}]^T$ denotes the coefficient vector. Therefore, we also apply the sparse regression method to learn the diffusion term.

Substituting the kernel function (\ref{kernel}) of the L\'evy jump measure into the expression for $S_{ij}^{\alpha ,\beta, \varepsilon}(z)$ from Corollary \ref{cor2}'s third formula yields
\begin{equation}\label{Sz1}
	S_{ij}^{\alpha ,\beta, \varepsilon}(z)= \left\{
	\begin{array}{ll}
		\sigma_i^{\alpha_i}(z) k_{\alpha_i} \beta_i \frac{\varepsilon ^{2-\alpha_i}} {2-\alpha_i}, \ \ & i=j, \\
		0, & i \neq j.
	\end{array} \right.
\end{equation}
Therefore, $S_{ij}^{\alpha ,\beta, \varepsilon}(z)$ can be numerically computed using the estimated parameters $\tilde{\alpha}_i$, $\tilde{\beta}_i$ and the learned noise intensity function $\sigma_{i,\theta}(x_i)$ from Subsection \ref{learnlevysubsec}.

Using the numerical approximation $\mathbbm{P} \left\{ \left. x(h)-z \in \Gamma  \right| x(0)=z \right\} \approx \frac{\hat{M}}{M}$ and assuming system ergodicity, Corollary \ref{cor2}'s third formula yields the following system of linear equations
\begin{equation}\label{eqdiffusion}
\begin{aligned}
  & A{d_{ij}}={{B}_{ij}}, \\
 & A=\left[ \begin{array}{ccc}
 	{{\psi }_{1}}\left( \hat{z}^{(1)} \right) & \cdots  & {{\psi }_{K}}\left( \hat{z}^{(1)} \right)  \\
 	\vdots  & \ddots  & \vdots   \\
 	{{\psi }_{1}}\left( \hat{z}^{(\hat{M})} \right) & \cdots  & {{\psi }_{K}}\left( \hat{z}^{(\hat{M})} \right)  \\
\end{array} \right], \\
 & {{B}_{ij}}=\hat{M}{{M}^{-1}}{{h}^{-1}}{{\left[ \left( \hat{x}_i^{(1)}- \hat{z}_i^{(1)} \right)\left( \hat{x}_j^{(1)}- \hat{z}_j^{(1)} \right),\ ...,\ \left( \hat{x}_i^{(\hat{M})}- \hat{z}_i^{(\hat{M})} \right)\left( \hat{x}_j^{(\hat{M})}- \hat{z}_j^{(\hat{M})} \right) \right]}^{T}} \\
 & \ \ \ \ \ \ \ - \left[ S_{ij}^{\tilde{\alpha}_i, \tilde{\beta}_i, \varepsilon} (\hat{z}^{(1)}), ..., S_{ij}^{\tilde{\alpha}_i, \tilde{\beta}_i, \varepsilon} (\hat{z}^{(\hat{M})}) \right]^T. \\
\end{aligned}
\end{equation}
This system can be transformed into a least squares regression problem
\begin{equation}\label{diffusionloss1}
	\tilde{d}_{ij}= \argmin_{d_{ij}} \| B_{ij}- A d_{ij} \|_2^2.
\end{equation}
To promote sparsity, we add a regularization term to Eq. (\ref{diffusionloss1}), yielding
\begin{equation}\label{diffusionloss2}
	\tilde{d}_{ij}= \argmin_{d_{ij}} \| B_{ij}- A d_{ij} \|_2^2 +\lambda \| d_{ij} \|_1.
\end{equation}

Similarly, we employ the combination algorithm of binning, stepwise sparse regression, and cross-validation for low-dimensional systems and SINDy for high-dimensional systems to learn the diffusion term. After processing all $n$ dimensions for $i$, we successfully learn the diffusion matrix from data. Due to the symmetry of the diffusion matrix, we only require $n(n+1)/2$ algorithm executions, iterating over $j=i, i+1, \ldots , n$ for each $i=1, 2, \ldots , n$.

\subsection{Convergence and error analysis}
\label{errorsubsec}
In this subsection, we perform convergence proof and error analysis for the previous proposed algorithms. For simplification of representation, we consider the system dimension $n=1$ so that we can omit the subscript $i$. For multi-dimensional systems, we can consider each dimension separately. We first make some assumptions.

\newtheorem{ass}{\bf Assumption} 
\begin{ass}\label{ass1}
	The sample size $M \to \infty$, the time step $h \to 0$, the binning number $N_c \to \infty$, and $MhN_c^{-1} \to \infty$.
\end{ass}

\begin{ass}\label{ass2}
	Truncation parameter $\varepsilon$ is fixed and $\varepsilon \gg h^{1/2}$.
\end{ass}

\begin{ass}\label{ass3}
	There exists a constant $\delta_0$ satisfying $0< \delta_0 <0.5$ such that $\max\{\alpha, |\alpha-1|, 2-\alpha\} >\delta_0$.
\end{ass}

\begin{ass}\label{ass4}
	The initial point dataset $\mathbbm{Z}$ are uniformly distributed, i.e., $M_l N_c M^{-1} \approx 1$ for $l=1,2,...,N_c$.
\end{ass}

Under these assumptions, we can derive the following theorem to establish the error rates for L\'evy parameter estimation.

\begin{thm}(Error rates for L\'evy parameter estimation)\\
	\label{thmerror1}
	Under the Assumptions \ref{ass1}-\ref{ass4}, we have the following error estimation:
	\begin{enumerate}[1)]
		\item Stability index:
		\begin{equation}\label{eqerror1}
		| \tilde{\alpha}- \alpha | = O(h) + O(N_c^{-1}) + O_{\mathbbm{P}} (\sqrt{M^{-1}h^{-1}N_c}).
		\end{equation}
		
		\item Skewness parameter:
		\begin{equation}\label{eqerror2}
			| \tilde{\beta}- \beta |= O(h) + O(N_c^{-1}) + O_{\mathbbm{P}} (\sqrt{M^{-1}h^{-1}N_c}).
		\end{equation}
		
		\item Noise intensity function:
		\begin{equation}\label{eqerror3}
			\frac{1}{N_c} \sum_{l=1}^{Nc} | \tilde{\sigma}(\hat{z}_l)- \sigma(\hat{z}_l) | = O(h) + O(N_c^{-1}) + O_{\mathbbm{P}} (\sqrt{M^{-1}h^{-1}N_c}).
		\end{equation}
	\end{enumerate}
\end{thm}

Here the notation '$O_{\mathbbm{P}}$' denotes stochastic boundedness. For the noise intensity function, the subsequent regression operation will introduce an additional regression error. Proof of Theorem \ref{thmerror1} are presented in Appendix \ref{App:H}. The following theorem establishes the convergence of learning algorithms for the drift and diffusion terms.

\begin{thm}(Convergence of drift and diffusion estimation)\\
	\label{thmerror2}
	Under the Assumptions \ref{ass1}-\ref{ass4} and the condition that all L\'evy parameters and functions are accurately identified, the solutions $\{ \tilde{c}, \tilde{d} \} \in \mathbbm{R}^{K_2}$ of the regression problems
	\begin{align}\label{converg1}
		& \tilde{c} = \argmin_{c \in \mathbbm{R}^{K_2}} \| B_b- A c \|_2^2, \\
		& \tilde{d} = \argmin_{d \in \mathbbm{R}^{K_2}} \| B_a- A d \|_2^2,
	\end{align}
    converge to the coefficient vectors of the best approximation problems
    \begin{align}\label{converg2}
    	& \tilde{c} = \argmin_{c \in \mathbbm{R}^{K_2}} \| b- \sum_{k=1}^{K_2} c_k\psi_k \|_{L_{\mu}^2}^2, \\
    	& \tilde{d} = \argmin_{d \in \mathbbm{R}^{K_2}} \| a- \sum_{k=1}^{K_2} d_k\psi_k \|_{L_{\mu}^2}^2,
    \end{align}
    in the space $L_{\mu}^2$ of square-integrable functions with respect to the measure $\mu$, where the subscripts '$b$' and '$a$' of $B$ indicate drift and diffusion.
\end{thm}

Proof of Theorem \ref{thmerror2} is inspired by Ref. \cite{Boninsegna2018} and appears in Appendix \ref{App:I}.

There are several remarks concerning these two theoretical results. First, note that the parameter $\alpha$ cannot be close to $\{ 0, 1, 2 \}$ (e.g., $\alpha=0.001$, $1.001$, or $1.999$) according to Assumption \ref{ass3}, since Eqs. (\ref{eq4}), (\ref{Rz1}) (i.e., the expression of $R^{\alpha, \beta, \varepsilon}$), and (\ref{Sz1}) (i.e., the expression of $S^{\alpha, \beta, \varepsilon}$) contain the factors $\alpha^{-1}$, $(1-\alpha)^{-1}$, and $(2-\alpha)^{-1}$, respectively. Specifically, if $\alpha$ is close to 0, the estimation of L\'evy noise parameters will be inaccurate; if $\alpha$ is close to 1 (or 2), the factor of $(1-\alpha)^{-1}$ (or $(2-\alpha)^{-1}$) will significantly amplify the $\sigma$ estimation error in $R^{\alpha, \beta, \varepsilon}$ (or $S^{\alpha, \beta, \varepsilon}$), which introduces an additional stochastic constant term into the inferred drift (or diffusion) term, causing constant errors. Numerical experiments show that with a data size of $10^7$, the error is acceptable as long as the distance between the parameter ? and $\{ 0, 1, 2 \}$  exceeds 0.05.

Second, the approximation of the drift and diffusion terms involves multiple sources of error. The first is the error propagated through $R^{\alpha, \beta, \varepsilon}$ or $S^{\alpha, \beta, \varepsilon}$ due to the preliminary error from the estimation of the L\'evy noise parameters. The second is the random fluctuations caused by estimating conditional expectations and probabilities with finite samples. The third is the bias introduced by ignoring the initial state in the conditional probabilities. The fourth is the bias resulting from using a finite set of basis functions to approximate an infinite-dimensional function space, as well as overfitting or underfitting.

Third, according to Assumption \ref{ass1} and Theorem \ref{thmerror1}, the error will converge to zero only if $M \to \infty$, $h \to 0$, $N_c \to \infty$, and $MhN_c^{-1} \to \infty$ are simultaneously satisfied. If $M$ and $h$ are held constant while only $N_c$ is increased, the algorithm's accuracy may even deteriorate. This is intuitively easy to understand because increasing $N_c$ leads to a reduction in the amount of data within each bin. As mentioned earlier, we can reduce the error caused by regional discretization by repeatedly selecting different values of $N_c$ when other parameters are fixed.

Fourth, under limiting conditions, we require $\varepsilon$ to be fixed. Numerically, $\varepsilon$ cannot be chosen too large (as it would cause large jumps to become too rare) or too small (as it would cause the estimation of L\'evy parameters to be affected by the diffusion term). Through extensive subsequent numerical experiments, it is found that the suitable range for $\varepsilon$ is 10 to 50 times $h^{1/2}$. Furthermore, the value of $\varepsilon$ used in the L\'evy parameter estimation algorithm does not necessarily need to be the same as that used in the drift and diffusion term estimation algorithms.

\section{Examples}
\label{exam}
In this section, we present two prototypical examples to illustrate our method for discovering stochastic dynamical systems with multiplicative (Gaussian) Brownian noise and multiplicative (non-Gaussian) L\'evy noise from sample path data.

\newtheorem{exa}{\bf Example}
\begin{exa}\label{exa1}
Consider the two-dimensional stochastic Maier-Stein system \cite{maier1993effect}
\begin{equation}\label{MSsystem}
\begin{aligned}
  & {\rm d}{{x}_{1}}=(x_1- x_1^3- 5x_1 x_2^2){\rm d}t+\left( 1+{{x}_1} \right) {\rm d}{{B}_{1,t}}+ {\rm d}{{B}_{2,t}}+ (1-x_1+x_1^2){\rm d}{{L}_{1,t}}, \\
  & {\rm d}{{x}_{2}}=-\left( 1+x_1^2 \right) x_2 {\rm d}t+{{x}_{2}}{\rm d} {{B}_{2,t}}+ \frac{1}{1+0.5x_2^2} {\rm d}{{L}_{2,t}}. \\
\end{aligned}
\end{equation}
The drift coefficient, the diffusion coefficient, and the L\'evy noise intensity functions to be identified take the following form
$$
\begin{aligned}
  & b\left( x \right)={{\left[ x_1- x_1^3- 5x_1 x_2^2,\ -\left( 1+x_1^2 \right) x_2 \right]}^{T}}, \\
  & a\left( x \right)=\left[ \begin{array}{cccc}
   2+2{{x}_{1}}+x_{1}^{2} & \ &{{x}_{2}} \\
   {{x}_{2}} & \ & x_{2}^{2} \\
\end{array} \right]. \\
  & \sigma_1(x)= 1-x_1+x_1^2, \\
  & \sigma_2(x)= \frac{1}{1+0.5x_2^2}. \\
\end{aligned}
$$
The noise parameters in L\'evy jump measure are chosen as ${{\alpha }_{1}}=0.5$, ${{\beta }_{1}}=0.5$, and ${{\alpha }_{2}}=1.5$, ${{\beta }_{2}}=-0.5$.

To simulate the system's sample data, we fix the time step $h=0.001$ and randomly select $M=10^7$ initial points $\mathbbm{Z}=\left[ {z_{1}}, {z_{2}}, ..., {z_{M}} \right]$ uniformly distributed in the domain $\left[ -2,2 \right]\times \left[ -2, 2 \right]$. The dataset $\mathbbm{X}$ is then generated by evolving the system from $\mathbbm{Z}$ for time $h$ using the Euler scheme.

Using the algorithm in Section \ref{learnlevysubsec}, we learn the L\'evy jump measure and noise intensity with fixed parameters $N=1$, $m=5$, and $\varepsilon=0.5$. The interval $[-2,2]$ in each dimension is divided into $N_c=10,11,...,25$ bins, yielding 280 values of the noise intensity function. From Eqs. (\ref{eq6}) and (\ref{eq7}), we compute two parameter sets ($\alpha$, $\beta$) and present the results in Table \ref{tab:1}. The learned results show excellent agreement with the true parameters.

\begin{table}[htbp]
  \centering
  \caption{Learned L\'evy noise parameters for the Maier-Stein system}
    \begin{tabular}{ccccc}
    \toprule
    \multicolumn{1}{c}{\multirow{2}[2]{*}{Parameter}} & \multicolumn{2}{c}{${L}_{1}$} & \multicolumn{2}{c}{${L}_{2}$}\\
          & True  & \multicolumn{1}{p{4.04em}}{Learned} & True  & \multicolumn{1}{p{4.04em}}{Learned} \\
    \midrule
    $\alpha$      & 0.5   & 0.5138 & 1.5   & 1.5182 \\
    $\beta$      & 0.5   & 0.5086 & -0.5  & -0.4785 \\
    \bottomrule
    \end{tabular}%
  \label{tab:1}%
\end{table}%

Using Eq. (\ref{eq8}), we approximate the L\'evy noise intensity function values (denoted by blue stars in Fig. \ref{fig1}(a,b)), which show good agreement with the true functions. For subsequent estimation of drift and diffusion terms, we learn the noise intensity functions from these approximations. First, we solve the sparse regression problem (\ref{sigmaloss2}) via stepwise sparse regressor with cross-validation, using the candidate function dictionary $\Phi(\xi) = \{ 1, \xi, \xi^2, \xi^3, \xi^4 \}$. The sparse coefficients are listed in Table \ref{tab:2}, with corresponding function curves shown in Fig. \ref{fig1}. Second, we approximate the functions using fully connected neural networks, with results also displayed in Fig. \ref{fig1}(a,b).

\begin{figure}
	\centering
	\includegraphics[width=12cm]{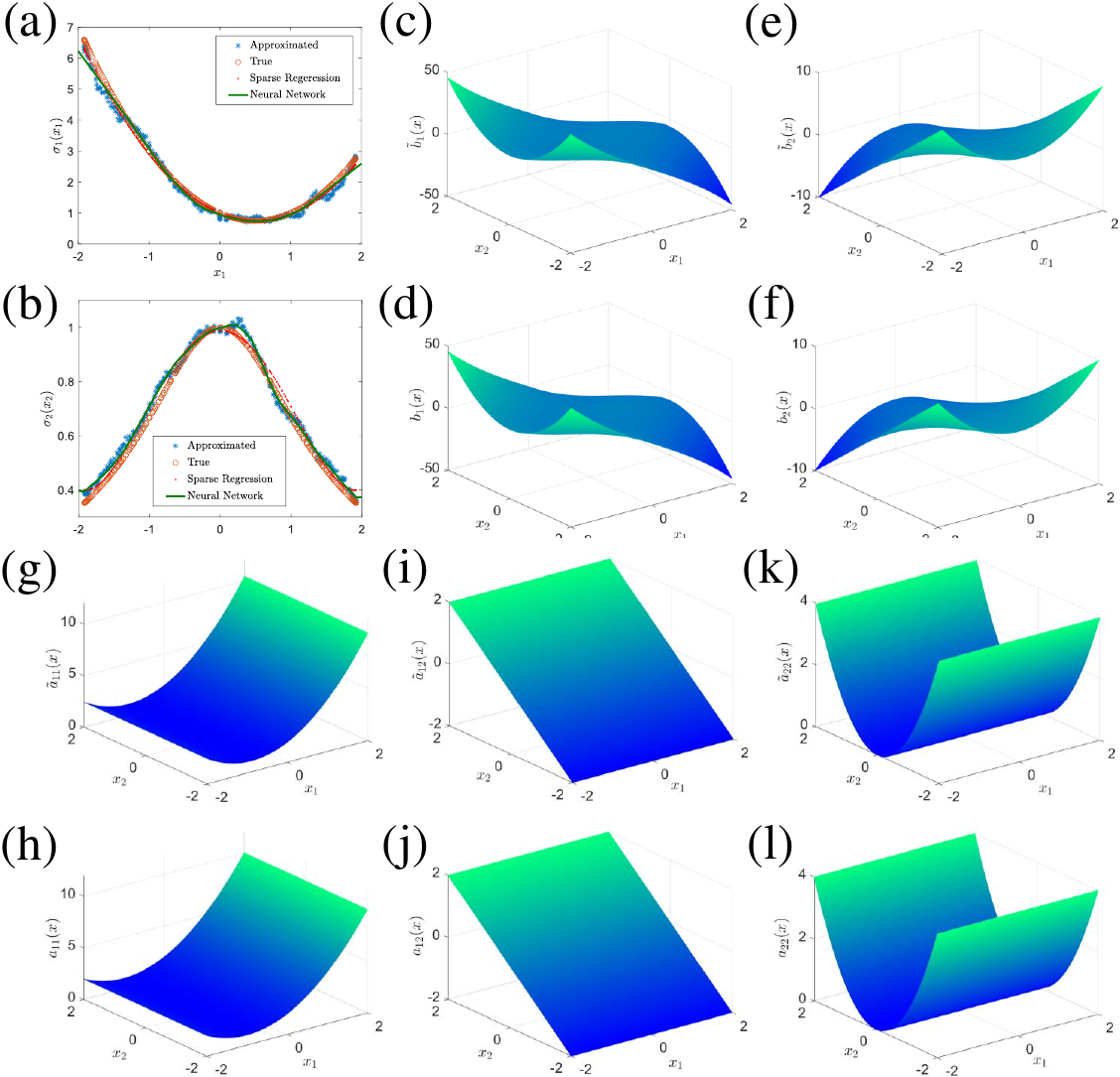}
	\caption{Comparison between learned and true coefficient functions of the Maier-Stein system. (a,b) The approximated values, sparse regression results, neural network results, and true functions of the L\'evy noise intensity functions $\sigma_1(x_1)$ and $\sigma_2(x_2)$, denoted by blue stars, red points, green curves, and orange circles, respectively. (c,e) The learned drift coeffcient $\tilde{b}(x)$. (d,f) The true drift coeffcient $b(x)$. (g,i,k) The learned diffusion coeffcient $\tilde{a}(x)$. (h,j,l) The true diffusion coeffcient $a(x)$.}
	\label{fig1}
\end{figure}

\begin{table}[htbp]
	\centering
	\caption{Learned L\'evy noise intensity functions for the Maier-Stein system}
	\begin{tabular}{cccc}
		\toprule
		\multicolumn{1}{c}{\multirow{2}[1]{*}{Basis}} & \multicolumn{2}{c}{${{\sigma}_{1}}\left( x_1 \right)$} & \multicolumn{1}{c}{${{\sigma}_{2}}\left( x_2 \right)$} \\
		& True  & \multicolumn{1}{p{4.04em}}{Learned} &  \multicolumn{1}{p{4.04em}}{Learned} \\
		\midrule
		\multicolumn{1}{c}{1}            & 1   & 1.0021   & 0.9875    \\
		\multicolumn{1}{c}{$\xi$}        & -1  & -0.9735  & 0         \\
		\multicolumn{1}{c}{$\xi^{2}$}    & 1   & 0.9312   & -0.3247   \\
		\multicolumn{1}{c}{$\xi^3$}      & 0   & 0        & 0         \\
		\multicolumn{1}{c}{$\xi^4$}      & 0   & 0        & 0.0448    \\
		\bottomrule
	\end{tabular}%
	\label{tab:2}%
\end{table}%

As shown in Fig. \ref{fig1}(a,b), both methods yield learned noise intensity functions that agree well with the true functions. Table \ref{tab:2} further demonstrates that the sparse coefficients learned for $\sigma_1(x_1)$ match the true values. Although $\sigma_2(x_2)$ is not a polynomial function, the sparse regression results obtained using polynomial basis functions remain reasonably accurate.

We define the approximation errors for sparse regression and neural networks as
\begin{equation}\label{errorsrnn}
	e_i^{sr}= \frac{\sqrt{\mathbbm{E} (\sigma_i^{sr} (x_i)- \sigma_i (x_i))^2}} {\max_{x_i} |\sigma_i (x_i)|}, \ e_i^{nn}= \frac{\sqrt{\mathbbm{E} (\sigma_i^{nn} (x_i)- \sigma_i (x_i))^2}} {\max_{x_i} |\sigma_i (x_i)|}, \ i=1,2.
\end{equation}
Using these definitions, we compute $e_1^{sr}=0.0179$, $e_2^{sr}=0.0287$, $e_1^{nn}=0.0255$, $e_2^{nn}=0.0282$. Thus, sparse regression achieves higher accuracy for polynomial functions ($i=1$), while neural networks perform better for non-polynomial cases ($i=2$).

\begin{table}[htbp]
  \centering
  \caption{Learned drift term for the Maier-Stein system}
    \begin{tabular}{ccccc}
    \toprule
    \multicolumn{1}{c}{\multirow{2}[1]{*}{Basis}} & \multicolumn{2}{c}{${{b}_{1}}\left( x \right)$} & \multicolumn{2}{c}{${{b}_{2}}\left( x \right)$} \\
          & True  & \multicolumn{1}{p{4.04em}}{Learned} & True  & \multicolumn{1}{p{4.04em}}{Learned} \\
    \midrule
    \multicolumn{1}{c}{1}                & 0   & 0       & 0     & 0    \\
    \multicolumn{1}{c}{${x}_{1}$}        & 1   & 1.0478  & 0     & 0 \\
    \multicolumn{1}{c}{${x}_{2}$}        & 0   & 0       & -1    & -1.0031 \\
    \multicolumn{1}{c}{$x_{1}^{2}$}      & 0   & 0       & 0     & 0   \\
    \multicolumn{1}{c}{${x}_{1}{x}_{2}$} & 0   & 0       & 0     & 0   \\
    \multicolumn{1}{c}{$x_{2}^{2}$}      & 0   & 0       & 0     & 0    \\
    \multicolumn{1}{c}{$x_1^3$}          & -1  & -1.0086 & 0     & 0   \\
    \multicolumn{1}{c}{$x_{1}^{2}x_2$}   & 0   & 0       & -1    & -1.0043 \\
    \multicolumn{1}{c}{$x_{1}x_2^{2}$}   & -5  & -4.9929 & 0     & 0   \\
    \multicolumn{1}{c}{$x_2^3$}          & 0   & 0       & 0     & 0   \\
    \bottomrule
    \end{tabular}%
  \label{tab:3}%
\end{table}%

We next learn the drift and diffusion terms using binning, stepwise sparse regressor, and cross-validation with the candidate dictionary $\Psi(x)= [1, x_1, x_2, x_1^2, x_1x_2, x_2^2, x_1^3, x_1^2x_2, x_1x_2^2, x_2^3]$ and parameter $\varepsilon=1$. Tables \ref{tab:3} and \ref{tab:4} present the resulting sparse coefficients for the drift and diffusion terms, respectively. Figs. \ref{fig1}(c-l) display the comparison between learned and true drift and diffusion functions. It is seen that all coefficients closely match the true values and the corresponding functional images agree well, confirming the method's effectiveness and accuracy.

\begin{table}[htbp]
  \centering
  \caption{Learned diffusion term for the Maier-Stein system}
    \begin{tabular}{ccccccc}
    \toprule
    \multicolumn{1}{c}{\multirow{2}[1]{*}{Basis}} & \multicolumn{2}{c}{${{a}_{11}}\left( x \right)$} & \multicolumn{2}{c}{${{a}_{12}}\left( x \right)$} & \multicolumn{2}{c}{${{a}_{22}}\left( x \right)$}  \\
     & True  & \multicolumn{1}{p{4.04em}}{Learned} & True  & \multicolumn{1}{p{4.04em}}{Learned} & True  & \multicolumn{1}{p{4.04em}}{Learned} \\
    \midrule
    \multicolumn{1}{c}{1}                  & 2     & 1.9682  & 0     & 0      & 0     & 0  \\
    \multicolumn{1}{c}{${x}_{1}$}          & 2     & 1.9892  & 0     & 0      & 0     & 0   \\
    \multicolumn{1}{c}{${x}_{2}$}          & 0     & 0       & 1     & 0.9954 & 0     & 0    \\
    \multicolumn{1}{c}{$x_{1}^{2}$}        & 1     & 1.1169  & 0     & 0      & 0     & 0    \\
    \multicolumn{1}{c}{${x}_{1}{x}_{2}$}   & 0     & 0       & 0     & 0      & 0     & 0    \\
    \multicolumn{1}{c}{$x_{2}^{2}$}        & 0     & 0       & 0     & 0      & 1     & 0.9894  \\
    \multicolumn{1}{c}{${x}_{1}^3$}        & 0     & 0       & 0     & 0      & 0     & 0    \\
    \multicolumn{1}{c}{$x_{1}^{2}x_2$}     & 0     & 0       & 0     & 0      & 0     & 0    \\
    \multicolumn{1}{c}{$x_{1}x_2^{2}$}     & 0     & 0       & 0     & 0      & 0     & 0    \\
    \multicolumn{1}{c}{$x_{2}^{3}$}        & 0     & 0       & 0     & 0      & 0     & 0    \\
    \bottomrule
    \end{tabular}%
  \label{tab:4}%
\end{table}%

Finally, we examine how data volume $M$ and time step $h$ affects algorithmic accuracy with $N_c$ and $\varepsilon$ fixed to verify the results of error analysis. Define the approximation errors for L\'evy noise parameters and intensity functions as
\begin{equation}\label{errorabs}
	e_i^{\alpha}= |\tilde{\alpha}_i-\alpha_i|, \ e_i^{\beta}= |\tilde{\beta}_i- \beta_i|, \ e_i^{\sigma}= \frac{\frac{1}{N_c} \sum_{l=1}^{Nc} |\tilde{\sigma}_i (\hat{z}_{i,l})- {{\sigma}_{i}} (\hat{z}_{i,l})|} {\max_{l} |{{\sigma}_{i}} (\hat{z}_{i,l}|} , \ i=1,2.
\end{equation}
The effects of the hyperparameters on these errors are illustrated in Fig. \ref{fig2}.

First, Fig. \ref{fig2}(a-c) show these errors across varying data volumes. All errors exhibit a decreasing trend with increasing data volume, confirming the algorithm's convergence. Note that the slopes of all errors decreasing with $M$ are consistent with $M^{-0.5}$, which aligns with the error analysis of Theorem \ref{thmerror1}, i.e., $O_{\mathbb{P}} (\sqrt{M^{-1}h^{-1}N_c})$. When $M > 10^8$, the rate of error reduction slows down, as the first two error terms $O(h) + O(N_c^{-1})$ dominate the errors now. Notably, under limited data conditions, $\sigma_2$ delivers superior accuracy to $\sigma_1$--a phenomenon originating from variations in $\alpha$: larger jump amplitudes with sparser occurrences due to smaller $\alpha$ amplify errors when data are insufficient. Additionally, by comparing the three figures, we observe that the errors for $\beta$ are the smallest, which stems from two reasons: on one hand, in the derivation of the error for $\beta$, error cancellation effects occur on both sides; on the other hand, the computation of $\beta$ involves a more comprehensive utilization of the data.

\begin{figure}
	\centering
	\includegraphics[width=12cm]{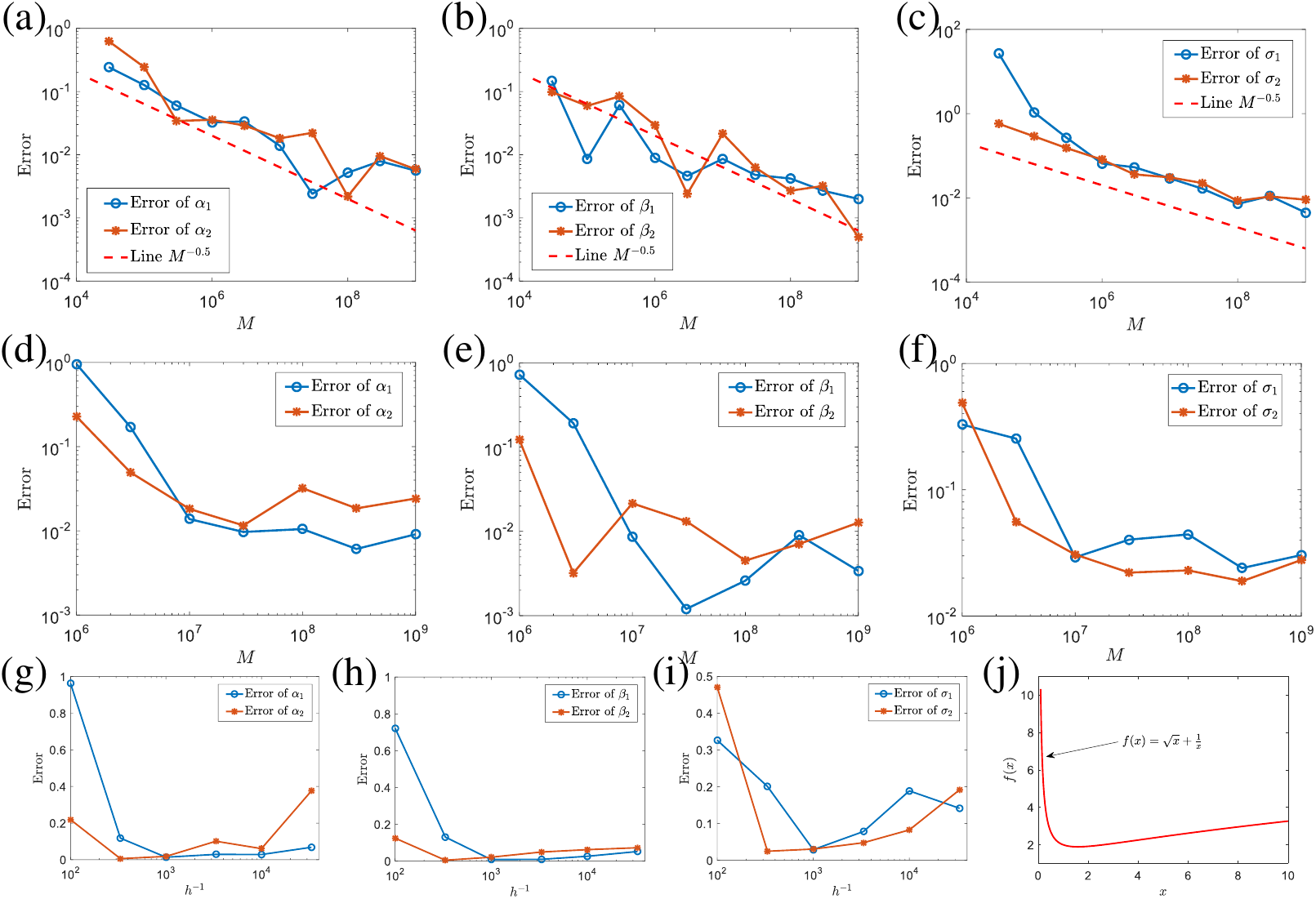}
	\caption{The influence of hyperparameters such as data volume $M$ and time step $h$ on approximation errors of L\'evy noise parameters $\alpha$, $\beta$, and $\sigma$, with $N_c$ and $\varepsilon$ fixed. The blue line with circles and the brown line with stars represent the first and second dimensions, respectively. (a-c) The approximation errors of L\'evy noise parameters about data volume $M$ with $h=0.001$ fixed. The red dashed lines indicate the function $20 \times M^{-0.5}$. (d-f) The approximation errors of L\'evy noise parameters about data volume $M$ with $Mh=10^4$ fixed. (g-i) The approximation errors of L\'evy noise parameters about $h^{-1}$ with $M=10^7$ fixed. (j) The graph of the function $f(x) = \sqrt{x} + \frac{1}{x}$.}
	\label{fig2}
\end{figure}

Second, Fig. \ref{fig2}(d-f) demonstrate the relationship between errors and the data volume $M$ under the condition of a fixed $Mh = 10^4$. It can be observed that when $M \geq 10^7$, the errors remain almost constant as $M$ increases, which further validates the error term $O_{\mathbb{P}} (\sqrt{M^{-1}h^{-1}N_c})$. For $M < 10^7$, the errors increase due to the limited data volume and the accompanying increase in $O(h)$.

Third, Fig. \ref{fig2}(g-i) illustrate the relationship between errors and $h^{-1}$ under the condition of a fixed $M = 10^7$. It is observed that the error curves initially decrease rapidly and then increase slowly, a pattern that aligns almost exactly with the shape of the function $f(x) = \sqrt{x} + \frac{1}{x}$ shown in Fig. \ref{fig2}(j). Therefore, the error exhibits a functional relationship with the time step $h$, given by $e(h) = c_1 \sqrt{h^{-1}} + c_2 (h^{-1})^{-1} = c_1 \sqrt{h^{-1}} + c_2 h$. This confirms the results of our error analysis $O(h) + O(N_c^{-1}) + O_{\mathbbm{P}} (\sqrt{M^{-1}h^{-1}N_c})$.
\end{exa}

\begin{exa}\label{exa2}
To illustrate our framework's capability in inferring high-dimensional non-Gaussian stochastic dynamical systems, we consider stochastic R\"ossler oscillator dynamics on the weighted network \cite{gao2024learning}
\begin{equation}\label{Rossler}
	\begin{aligned}
		& {\rm d}{{x}_{i,1}}= (-x_{i,2}- x_{i,3}+ \sum_{j=1}^J F_{i,j} (x_{j,1}-x_{i,1})) {\rm d}t+ x_{i,1} {\rm d}{{B}_{i,1,t}}+ (1+0.5x_{i,1}^2) {\rm d}{{L}_{i,1,t}}, \\
		& {\rm d}{{x}_{i,2}}= \left( x_{i,1}+ 0.35x_{i,2} \right) {\rm d}t+ {\sqrt{1+x_{i,2}^2}} {\rm d} {{B}_{i,2,t}}+ \frac{1} {{1+0.5x_{i,2}^2}} {\rm d}{{L}_{i,2,t}}, \\
		& {\rm d}{{x}_{i,3}}= \left( 0.2+ x_{i,1}x_{i,3}- 5.7x_{i,3} \right) {\rm d}t+ \sqrt{1+0.5x_{i,2}^2+x_{i,3}^2} {\rm d} {{B}_{i,3,t}}+ ({1+\sin^2{x_{i,3}}}) {\rm d}{{L}_{i,3,t}}. \\
	\end{aligned}
\end{equation}
Here, $J$ represents the number of nodes in the network, and the system dimensionality is $n=3J$. We set $J=5$, with the matrix $F$ defined as
$$
F= \left[ \begin{array}{ccccc}
		0 & 0.7161 & 0 & 0 & 0 \\
		0 & 0 & -1.2678 & 0 & 0 \\
		0 & 0 & 0 & -1.0141 & 0 \\
		0 & 0 & 0 & 0 & 2.3633 \\
		2.7307 & 0 & 0 & 0 & 0 \\
	\end{array} \right],
$$
where the nonzero elements are randomly generated. The noise parameters in L\'evy jump measure are specified as ${{\alpha }_{i,1}}=0.5$, ${{\beta }_{i,1}}=0.5$, ${{\alpha }_{i,2}}=1.1$, ${{\beta }_{i,2}}=0$, ${{\alpha }_{i,3}}=1.5$, ${{\beta }_{i,3}}=-0.5$ for $i=1,2,3,4,5$.

We also set the time step to $h=0.001$ and select $M_0=20000$ initial points $\mathbbm{Z}_0=\left[ {z_{1}}, {z_{2}}, ..., {z_{M_0}} \right]$ randomly and uniformly distributed in the domain $\left[ -2,2 \right] \times \left[ -2,2 \right] \times \cdots \times \left[ -2, 2 \right]$. Each initial point is integrated over the time interval $T=1$. The data sets $\mathbbm{Z}$ and $\mathbbm{X}$, generated via the Euler scheme, each contain $M=2 \times 10^7$ elements. To prevent divergence in simulated paths due to occasional large jumps in L\'evy motion, trajectories deviating from the study region are relocated to a randomly generated point within this domain.

We apply the algorithm in Section \ref{learnlevysubsec} to infer the L\'evy jump measure and noise intensity. Parameters are fixed at $N=1$, $m=5$, and $\varepsilon=0.5$. The interval $[-2,2]$ along each dimension is partitioned into $N_c=10,11,...,25$ bins, yielding 280 values of the noise intensity function. Using Eqs. (\ref{eq6}) and (\ref{eq7}), we compute 15 sets of parameters $\alpha$ and $\beta$, with results listed in Table \ref{tab:5}. All inferred values closely match the true parameters, and nearly all errors fall within $5\%$.

\begin{table}[htbp]
	\centering
	\caption{Learned L\'evy noise parameters for the R\"ossler system}
	\begin{tabular}{ccccccc}
		\toprule
		Parameter & $\alpha_{i,1}$ & $\alpha_{i,2}$ & $\alpha_{i,3}$ & $\beta_{i,1}$ & $\beta_{i,2}$ & $\beta_{i,3}$  \\
		\midrule
		True      & 0.5   & 1.1  & 1.5   & 0.5 & 0 & -0.5  \\
		Learned ($i=1$)   & 0.4877   & 1.1224  & 1.5547   & 0.5075 & 0.0093 & -0.5045  \\
		Learned ($i=2$)   & 0.5079   & 1.1250  & 1.5245   & 0.5005 & 0.0094 & -0.4930  \\
		Learned ($i=3$)   & 0.5017   & 1.0798  & 1.5376   & 0.4940 & 0.0055 & -0.4990  \\
		Learned ($i=4$)   & 0.4880   & 1.1234  & 1.5103   & 0.5133 & 0.0053 & -0.4890  \\
		Learned ($i=5$)   & 0.5034   & 1.1030  & 1.5372   & 0.4988 & 0.0124 & -0.4957  \\
		\bottomrule
	\end{tabular}%
	\label{tab:5}%
\end{table}%

Similar to the previous example, we estimate L\'evy noise intensity function values using Eq. (\ref{eq8}), denoted by blue stars in Fig. \ref{fig3}. We also employ two methods (stepwise sparse regression and neural networks) to learn these functions, with results shown in Fig. \ref{fig3}. The candidate function dictionary for sparse regression is chosen as $\Phi(\xi) = \{ 1, \xi, \xi^2, \xi^3, \xi^4 \}$. Table \ref{tab:6} lists the L\'evy noise intensity functions learned via sparse regression. As shown in Fig. \ref{fig3} and Table \ref{tab:6}, the inferred results match the true functions for both polynomial cases and complex forms such as rational or trigonometric functions.

\begin{table}[htbp]
	\centering
	\caption{Learned L\'evy noise intensity for the R\"ossler system}
	\begin{tabular}{cccc}
		\toprule
		Function & $\sigma_{i,1}(x)$ & $\sigma_{i,2}(x)$ & $\sigma_{i,3}(x)$  \\
		\midrule
		True      & $1+0.5x_{i,1}^2$   & $\frac{1} {{1+0.5x_{i,2}^2}}$  & ${1+\sin^2{x_{i,3}}}$  \\
		Learned ($i=1$)   & $0.9376+0.6279x_{1,1}^2$   & $0.9692-0.3270x_{1,2}^2+0.0497x_{1,2}^4$  & $1.0826+0.7745x_{1,3}^2-0.1401x_{1,3}^4$  \\
		Learned ($i=2$)   & $0.9994+0.4407x_{2,1}^2$   & $0.9886-0.3498x_{2,2}^2+0.0547x_{2,2}^4$  & $1.0612+0.8115x_{2,3}^2-0.1589x_{2,3}^4$   \\
		Learned ($i=3$)   & $1.0003+0.4824x_{3,1}^2$   & $0.9861-0.3472x_{3,2}^2+0.0501x_{3,2}^4$  & $1.0644+0.8029x_{3,3}^2-0.1556x_{3,3}^4$  \\
		Learned ($i=4$)   & $1.0538+0.4346x_{4,1}^2$   & $0.9804-0.3355x_{4,2}^2+0.0484x_{4,2}^4$  & $1.0403+0.8410x_{4,3}^2-0.1677x_{4,3}^4$  \\
		Learned ($i=5$)   & $0.9224+0.5447x_{5,1}^2$   & $0.9930-0.3472x_{5,2}^2+0.0491x_{5,2}^4$  & $1.0679+0.7910x_{5,3}^2-0.1519x_{5,3}^4$  \\
		\bottomrule
	\end{tabular}%
	\label{tab:6}%
\end{table}%

As shown in Fig. \ref{fig3}, the algorithm's approximations exhibit higher accuracy in the central region but slightly larger errors at the endpoints. This occurs because larger coordinate values at the extremities lead to increased drift terms, consequently affecting the estimation of L\'evy noise intensity. For dimensions $x_{i,1}$ and $x_{i,2}$, $i=1,2,3,4,5$, this error does not significantly affect the sparse regression results. However, for dimensions $x_{i,3}$, the endpoints coincide precisely with inflection points of the function, resulting in certain errors in sparse regression. This further leads to corresponding estimation errors for $a_{ii,3}(x)$ in Table \ref{tab:8}.

\begin{figure}
	\centering
    \includegraphics[width=12cm]{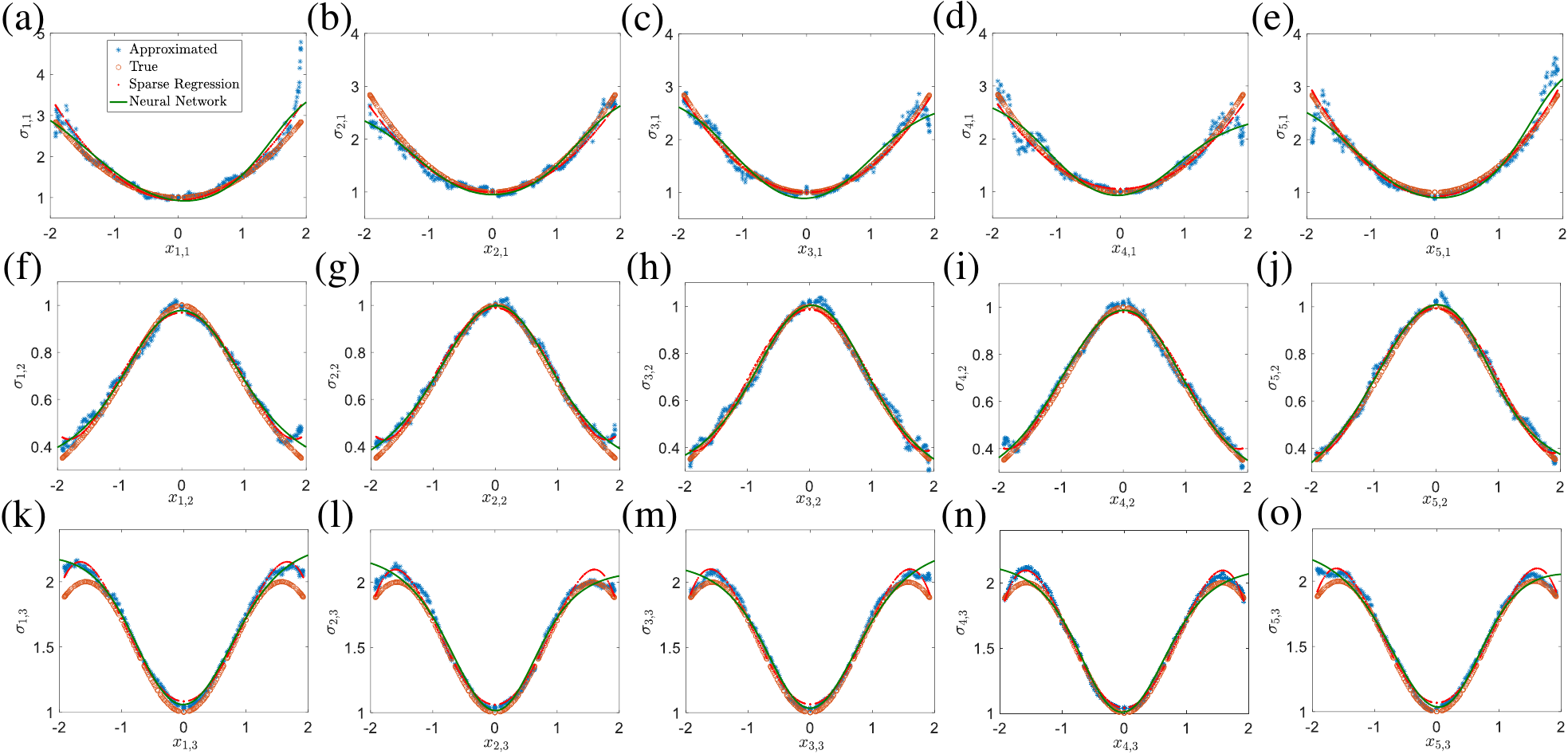}
	\caption{Comparison between the L\'evy noise intensity functions $\sigma_{i,j}(x_{i,j})$, $i=1,2,3,4,5$, $j=1,2,3$ of approximated values, sparse regression results, neural network results, and true functions of the R\"ossler system, denoted by blue stars, red points, green curves, and orange circles, respectively.}
	\label{fig3}
\end{figure}

We next employ sparse regression to learn the drift and diffusion terms using a second-order polynomial candidate dictionary with $\varepsilon=1$, resulting in a 136-element dictionary. The inferred drift and diffusion terms are presented in Tables \ref{tab:7} and \ref{tab:8}, respectively. Off-diagonal elements of the learned diffusion matrix are identically zero and thus omitted. All results closely approximate the true functions, demonstrating our method's effectiveness for the 15-dimensional R\"ossler system.

\begin{table}[htbp]
	\centering
	\caption{Learned drift term for the R\"ossler system}
	\begin{tabular}{cccc}
		\toprule
		Function & $b_{i,1}(x)$ & $b_{i,2}(x)$ & $b_{i,3}(x)$  \\
		\midrule
		True      & $-x_{i,2}- x_{i,3}+ \sum_{j=1}^5 F_{i,j} (x_{j,1}-x_{i,1})$   & $x_{i,1}+ 0.35x_{i,2}$  & $0.2+ x_{i,1}x_{i,3}- 5.7x_{i,3}$  \\
		Learned ($i=1$)   & $-0.9854x_{1,2}- 0.9941x_{1,3}+ 0.7003x_{2,1}- 0.7086x_{1,1}$   & $0.9809x_{1,1}+ 0.3521x_{1,2}$  & $0.2431+ 0.9921x_{1,1}x_{1,3}- 5.6306x_{1,3}$  \\
		Learned ($i=2$)   & $-0.9811x_{2,2}- 0.9819x_{2,3}- 1.2552x_{3,1}+ 1.2469x_{2,1}$   & $1.0059x_{2,1}+ 0.3540x_{2,2}$  & $0.2142+ 1.0257x_{2,1}x_{2,3}- 5.6240x_{2,3}$   \\
		Learned ($i=3$)   & $-0.9992x_{3,2}- 0.9737x_{3,3}- 1.0181x_{4,1}+ 1.0025x_{3,1}$   & $0.9770x_{3,1}+ 0.3476x_{3,2}$  & $0.2249+ 1.0113x_{3,1}x_{3,3}- 5.6187x_{3,3}$  \\
		Learned ($i=4$)   & $-0.9998x_{4,2}- 0.9932x_{4,3}+ 2.3295x_{5,1}- 2.3321x_{4,1}$   & $1.0013x_{4,1}+ 0.3460x_{4,2}$  & $0.2224+ 0.9593x_{4,1}x_{4,3}- 5.6304x_{4,3}$  \\
		Learned ($i=5$)   & $-0.9890x_{5,2}- 0.9861x_{5,3}+ 2.7161x_{1,1}- 2.7334x_{5,1}$   & $1.0365x_{5,1}+ 0.3427x_{5,2}$  & $0.2666+ 0.9741x_{5,1}x_{5,3}- 5.6450x_{5,3}$  \\
		\bottomrule
	\end{tabular}%
	\label{tab:7}%
\end{table}%

\begin{table}[htbp]
	\centering
	\caption{Learned diffusion term for the R\"ossler system}
	\begin{tabular}{cccc}
		\toprule
		Function & $a_{ii,1}(x)$ & $a_{ii,2}(x)$ & $a_{ii,3}(x)$  \\
		\midrule
		True      & $x_{i,1}^2$   & $1+x_{i,2}^2$  & $1+0.5x_{i,2}^2+x_{i,3}^2$  \\
		Learned ($i=1$)   & $0.9877x_{1,1}^2$   & $0.9629+0.9868x_{1,2}^2$  & $0.7902+0.5264x_{1,2}^2+0.7956x_{1,3}^2$  \\
		Learned ($i=2$)   & $0.9885x_{2,1}^2$   & $0.9547+0.9964x_{2,2}^2$  & $0.8786+0.5241x_{2,2}^2+0.8920x_{2,3}^2$   \\
		Learned ($i=3$)   & $0.9877x_{3,1}^2$   & $0.9869+0.9903x_{3,2}^2$  & $0.8393+0.5323x_{3,2}^2+0.8605x_{3,3}^2$  \\
		Learned ($i=4$)   & $0.9982x_{4,1}^2$   & $0.9571+0.9955x_{4,2}^2$  & $0.9293+0.5229x_{4,2}^2+0.9067x_{4,3}^2$  \\
		Learned ($i=5$)   & $0.9985x_{5,1}^2$   & $0.9713+0.9923x_{5,2}^2$  & $0.8444+0.5212x_{5,2}^2+0.8697x_{5,3}^2$  \\
		\bottomrule
	\end{tabular}%
	\label{tab:8}%
\end{table}%

Two remarks are noteworthy for this example. Firstly, based on the previous discussion, $\alpha$ should not be too close to the set $\{0,1,2\}$, so we set the second $\alpha$ parameter to 1.1 instead of 1. Secondly, while the 2-dimensional example uses $10^7$ data points, this 15-dimensional case employs $2 \times 10^7$. The proposed algorithms thus effectively overcome the curse of dimensionality, as data requirements increase minimally with dimensionality.
\end{exa}

\section{Discussion}
\label{disc}
This work presents a comprehensive framework for learning stochastic dynamical systems driven by both multiplicative (Gaussian) Brownian noise and multiplicative (non-Gaussian) L\'evy noise from sample path data. The core theoretical contribution lies in establishing the nonlocal Kramers-Moyal formulas (Theorem \ref{thm1}, Corollary \ref{cor2}, Theorem \ref{thm3}, Corollary \ref{cor4}), which generalize the classical Kramers-Moyal relations to accommodate the significant jumps characteristic of L\'evy motions. These formulas provide a rigorous foundation for connecting the short-time evolution of the system's transition probability density, or equivalently the statistics of sample paths, directly to the underlying SDE coefficients, namely the drift vector $b(x)$, diffusion matrix $a(x)$, L\'evy jump measure kernel $W(\xi)$, and noise intensity functions $\sigma(x)$.

Building upon this theory, we developed practical data-driven algorithms capable of simultaneously identifying all these governing components. The key features and implications of this work are:
\begin{itemize}
	\item Handling multiplicative and asymmetric L\'evy noise: The algorithms explicitly address the challenges posed by multiplicative noise intensity ($\sigma(x)$ dependent on state $x$) and asymmetric ($\beta \neq 0$) L\'evy noise, which are prevalent in complex systems but often neglected in simpler identification schemes.
	
	\item Convergence and error analysis: We prove the convergence results of the estimation algorithms for the drift and diffusion terms (Theorem \ref{thmerror2}), and establish an error analysis for the parameter estimation results under L\'evy noise (Theorem \ref{thmerror1}). Numerical experiments validate the theoretical findings regarding the order of errors with respect to the hyperparameters.
	
	\item Avoiding the curse of dimensionality (for component-wise $\sigma_i(x_i)$): By exploiting the component-wise structure of the L\'evy jump measure estimation (Corollary \ref{cor2}) and assuming noise intensity $\sigma_i(x)$ depends only on its corresponding coordinate $x_i$, the algorithm for learning the L\'evy measure scales efficiently with dimension. This is a significant advantage for moderately high-dimensional systems, as demonstrated by the successful identification in the 15-dimensional R\"ossler network (Example \ref{exa2}).
	
	\item Flexibility in function representation: The framework incorporates both sparse regression (using predefined dictionaries like polynomials) and neural networks for learning the functional forms of $b(x)$, $a(x)$, and $\sigma_i(x_i)$. This flexibility allows tailoring the approach to the complexity of the underlying system, as evidenced by the accurate reconstruction of diverse functions (polynomial, rational, trigonometric) in the examples (Tables \ref{tab:2}, \ref{tab:6}, Figs. \ref{fig1}, \ref{fig3}).
\end{itemize}

However, there are also some limitations in the proposed methods that need further considerations in the future:
\begin{itemize}
	\item General noise intensity dependence: The current efficiency for L\'evy measure identification relies on the assumption $\sigma_i(x)=\sigma_i(x_i)$. Relaxing this to allow $\sigma_i(x)$ to depend on all state variables necessitates discretizing the full state space, which suffers from the curse of dimensionality.
	
	\item Error Propagation: The algorithms involve sequential estimation: L\'evy parameters/intensity first, then drift/diffusion using these estimates. Errors in the first step propagate to the second. Investigating robust joint estimation schemes or error-correcting techniques could improve overall accuracy.
\end{itemize}

This work provides a principled and practical toolbox for discovering complex stochastic dynamics from data, particularly when driven by non-Gaussian fluctuations exhibiting large jumps and state-dependent intensities. The derived nonlocal Kramers-Moyal formulas establish a fundamental theoretical link between path properties and governing laws for stochastic systems with L\'evy noise. The developed algorithms open doors to data-driven modeling in diverse fields where such noise is prevalent, including:
\begin{itemize}
	\item Climate science: Modeling extreme weather events or abrupt climate shifts.
	
	\item Neuroscience: Analyzing neuronal spike trains or network dynamics with heavy-tailed inputs.
	
	\item Epidemiology: Capturing the impact of superspreading events in disease transmission models.
	
	\item Financial Mathematics: Modeling asset returns with heavy tails and volatility clustering.
	
	\item Biological Physics: Understanding molecular motor dynamics or gene expression bursting under non-Gaussian noise.
\end{itemize}
By enabling the extraction of interpretable SDEs with both Gaussian and non-Gaussian multiplicative noise from observations, this framework significantly advances our ability to understand, predict, and potentially control complex systems operating under the influence of discontinuous, heavy-tailed fluctuations.

\section*{Acknowledgement}
This research was supported by the National Natural Science Foundation of China (Grant Nos. 12302035, 12572038, and 12372030) and Fundamental Research Funds for the Central Universities (Grant No. 30925010411).  This work was also partly supported by the NSFC  International Collaboration Fund for Creative Research Teams (Grant W2541005), the Guangdong Provincial Key Laboratory of Mathematical and Neural Dynamical Systems (Grant 2024B1212010004),  the Cross Disciplinary Research Team on Data Science and Intelligent Medicine (2023KCXTD054),   the Guangdong-Dongguan Joint Research Fund (Grant 2023A151514 0016),  Dongguan Key Laboratory for Data Science and Intelligent Medicine, and the Guangdong-Dongguan Joint Research Grant 2023A1515140016.

\section*{Data Availability Statement}
The data that support the findings of this study are openly available in GitHub \url{https://github.com/liyangnuaa/NKM-for-SDE-with-multiplicative-Levy-noise}.

\appendix

 \section{$\alpha$-stable L\'evy processes}
\label{App:A}
A scalar $\alpha$-stable L\'evy process ${{L}_{t}}$ is a stochastic process with the following conditions:
\begin{enumerate}[i)]
\item ${{L}_{0}}=0$, a.s.;
\item Independent increments: for any choice of $n\ge 1$ and ${{t}_{0}}<{{t}_{1}}<\cdots <{{t}_{n-1}}<{{t}_{n}}$, the random variables ${{L}_{{{t}_{0}}}}$, ${{L}_{{{t}_{1}}}}-{{L}_{{{t}_{0}}}}$, ${{L}_{{{t}_{2}}}}-{{L}_{{{t}_{1}}}}$, $\cdots $, ${{L}_{{{t}_{n}}}}-{{L}_{{{t}_{n-1}}}}$ are independent;
\item Stationary increments: ${{L}_{t}}-{{L}_{s}}\sim {{S}_{\alpha }}\left( {{\left( t-s \right)}^{{1}/{\alpha }\;}},\beta ,0 \right)$;
\item Stochastically continuous sample paths: for every $s>0$, ${{L}_{t}}\to {{L}_{s}}$ in probability, as $t\to s$.
\end{enumerate}

The $\alpha$-stable L\'evy motion is a special yet widely used type of L\'evy process, defined by a stable random variable with distribution ${{S}_{\alpha }}\left( \delta ,\ \beta ,\ \lambda  \right)$ \cite{HeavytailBook, LMBook2, OEBbook}. Typically, $\alpha \in \left( 0,\ 2 \right]$ is the stability parameter, $\delta \in \left( 0,\ \infty  \right)$ the scaling parameter, $\beta \in \left[ -1,\ 1 \right]$ the skewness parameter and $\lambda \in \left( -\infty ,\ \infty  \right)$ the shift parameter.

For a stable random variable $X$ with $0<\alpha <2$, the heavy-tail behavior is characterized by
$$
\underset{x\to \infty }{\mathop{\lim }}\,{{y}^{\alpha }}\mathbbm{P}\left( X>y \right)={{C}_{\alpha }}\frac{1+\beta }{2}{{\delta }^{\alpha }};
$$
where ${{C}_{\alpha }}$ is a positive constant depending on $\alpha$. This implies polynomial decay of the tail probability. In $\alpha$-stable L\'evy motion, smaller $\alpha$ ($0<\alpha <1$) yields larger jumps with lower frequencies, while larger $\alpha$ ($1<\alpha <2$) produces smaller jumps with higher frequencies. The case $\alpha=2$ corresponds to (Gaussian) Brownian motion. Further details on L\'evy processes are found in Refs. \cite{Duan2015, applebaum2009levy}.

\section{Proof of Theorem \ref{thm1}}
\label{App:B}
1)  Since $p\left( x,0|z,0 \right)=\delta \left( x-z \right)$, then ${{p}_{i}}\left( {{x}_{i}},t|z,0 \right)=\delta \left( {{x}_{i}}-{{z}_{i}} \right)$ and ${{p}_{i}}\left( {{x}_{i}},0|z,0 \right)=0$ for arbitrary ${{x}_{i}}$ and  ${{z}_{i}}$ satisfying $\left| {{x}_{i}}-{{z}_{i}} \right|>\varepsilon $. As before, denote ${x^{i}}$ as the rest of the vector $x$ with  ${{x}_{i}}$ being removed. Then
$$
\begin{aligned}
  & \underset{t\to 0}{\mathop{\lim }}\,{{t}^{-1}}{{p}_{i}}\left( {{x}_{i}},t|z,0 \right) \\
  & =\underset{t\to 0}{\mathop{\lim }}\,{{t}^{-1}}\left[ {{p}_{i}}\left( {{x}_{i}},t|z,0 \right)-{{p}_{i}}\left( {{x}_{i}},0|z,0 \right) \right] \\
  & =\int_{{{\mathbbm{R}}^{n-1}}}{{{\left. \frac{\partial p\left( x,t|z,0 \right)}{\partial t} \right|}_{t=0}}{\rm d}{x^{i}}}\\
  & =-\int_{{{\mathbbm{R}}^{n-1}}}{\nabla \cdot \left[bp\left( x,0|z,0 \right) \right]{\rm d}{x^{i}}}+\frac{1}{2}\int_{{{\mathbbm{R}}^{n-1}}}{Tr\left[ H\left( ap\left( x,0|z,0 \right) \right) \right]{\rm d}{x^{i}}} \\
  & \ \  +\sum\limits_{j=1}^{n}{\int_{{{\mathbbm{R}}^{n-1}}}{\int_{\mathbbm{R}\backslash \left\{ 0 \right\}}{\left[p\left( x-{{\sigma }_{j}}(z){y}_{j}{e_{j}},0|z,0 \right)- p\left( x,0|z,0 \right)+{{\sigma }_{j}}(z)\chi _{j}^{\alpha }\left( {y}_{j} \right){y}_{j}\frac{\partial }{\partial {{x}_{j}}}p\left( x,0|z,0 \right) \right]W_{j}^{\alpha ,\beta }\left( {y}_{j} \right){\rm d}{y}_{j}} {\rm d}{x^{i}}}}.
\end{aligned}
$$
Since $\left| {{x}_{i}}-{{z}_{i}} \right|>\varepsilon $, $p\left( x,0|z,0 \right)=0$. Thus the above equation is reduced as
$$
\sum\limits_{j=1}^{n}{\int_{{{\mathbbm{R}}^{n-1}}}{\int_{\mathbbm{R}\backslash \left\{ 0 \right\}}{p\left( x-{{\sigma }_{j}}(z){y}_{j}{e_{j}},0|z,0 \right) W_{j}^{\alpha ,\beta }\left( {y}_{j} \right){\rm d}{y}_{j}} {\rm d} {x^{i}}}}.
$$
Since
$$
\begin{aligned}
  & p\left( x-{{\sigma }_{j}}(z){y}_{j}{e_{j}},0|z,0 \right) \\
  & =\sigma _{j}^{-1}(z)\delta \left( {{x}_{1}}-{{z}_{1}} \right)\cdots \delta \left( {{x}_{j-1}}-{{z}_{j-1}} \right)\delta \left( {{x}_{j}}-{{\sigma }_{j}}{y}_{j}-{{z}_{j}} \right)\delta \left( {{x}_{j+1}}-{{z}_{j+1}} \right)\cdots \delta \left( {{x}_{n}}-{{z}_{n}} \right), \\
\end{aligned}
$$
its integral is equal to zero for $j\ne i$ due to $\delta \left( {{x}_{i}}-{{z}_{i}} \right)=0$. Hence, we have
$$
\sum\limits_{j=1}^{n}{\int_{{{\mathbbm{R}}^{n-1}}}{\int_{\mathbbm{R}\backslash \left\{ 0 \right\}}{p\left( x-{{\sigma }_{j}}(z){y}_{j}{e_{j}},0|z,0 \right) W_{j}^{\alpha ,\beta }\left( {y}_{j} \right){\rm d}{y}_{j}} {\rm d} {x^{i}}}} =\sigma _{i}^{-1}(z)W_{i}^{\alpha ,\beta }\left( \sigma_{i}^{-1}(z) \left( {{x}_{i}}-{{z}_{i}} \right) \right).
$$

\noindent 2)  According to the Fokker-Planck equation (\ref{FPK}), we have
$$
\begin{aligned}
  & \underset{t\to 0}{\mathop{\lim }}\,{{t}^{-1}}\int_{x-z\in \Gamma }{\left( {{x}_{i}}-{{z}_{i}} \right)p\left( x,t|z,0 \right){\rm d}x} \\
  & =\underset{t\to 0}{\mathop{\lim }}\,{{t}^{-1}}\int_{x-z\in \Gamma }{\left( {{x}_{i}}-{{z}_{i}} \right)\left[ p\left( x,t|z,0 \right)-p\left( x,0|z,0 \right) +p\left( x,0|z,0 \right) \right]{\rm d}x} \\
  & =\int_{x-z\in \Gamma }{\left( {{x}_{i}}-{{z}_{i}} \right){{\left. \frac{ \partial p\left( x,t|z,0 \right)}{\partial t} \right|}_{t=0}} {\rm d}x} +\underset{t\to 0}{\mathop{\lim }}\,\left[ {{t}^{-1}}\int_{x-z\in \Gamma }{\left( {{x}_{i}}-{{z}_{i}} \right)\delta \left( x-z \right){\rm d}x} \right] \\
  & =-\sum\limits_{j=1}^{n}{\int_{x-z\in \Gamma }{\left( {{x}_{i}}-{{z}_{i}} \right)\frac{\partial }{\partial {{x}_{j}}}\left[ {{b}_{j}}\left( x \right) p\left( x,0|z,0 \right) \right]{\rm d}x}} \\
  & \ \ \ +\frac{1}{2}\sum\limits_{k,l=1}^{n}{\int_{x-z\in \Gamma }{\left( {{x}_{i}}-{{z}_{i}} \right)\frac{{{\partial }^{2}}}{\partial {{x}_{k}} \partial {{x}_{l}}}\left[ {{a}_{kl}}\left( x \right)p\left( x,0|z,0 \right) \right]{\rm d}x}} \\
  & \ \ \ +\sum\limits_{j=1}^{n}{\int_{x-z\in \Gamma }{\int_{\mathbbm{R} \backslash \left\{ 0 \right\}}{\left( {{x}_{i}}-{{z}_{i}} \right)\left[ p\left( x-{{\sigma }_{j}}(z){y}_{j}{e_{j}},0|z,0 \right)- p\left( x,0|z,0 \right)+{{\sigma }_{j}}(z)\chi _{j}^{\alpha }\left( {y}_{j} \right){y}_{j} \frac{\partial }{\partial {{x}_{j}}}p\left( x,0|z,0 \right) \right] W_{j} ^{\alpha ,\beta }\left( {y}_{j} \right){\rm d}{y}_{j}}{\rm d}x}}. \\
\end{aligned}
$$
The application of integration by parts into the first term leads to
$$
\begin{aligned}
  & \sum\limits_{j=1}^{n}{\int_{x-z\in \Gamma }{\left( {{x}_{i}}-{{z}_{i}} \right)\frac{\partial }{\partial {{x}_{j}}}\left[ {{b}_{j}}\left( x \right) p\left( x,0|z,0 \right) \right]{\rm d}x}} \\
  & =-\sum\limits_{j=1}^{n}{\int_{x-z\in \Gamma }{{{b}_{j}}\left( x \right) p\left( x,0|z,0 \right)\frac{\partial }{\partial {{x}_{j}}}\left( {{x}_{i}} -{{z}_{i}} \right){\rm d}x}} \\
  & =-\sum\limits_{j=1}^{n}{\int_{x-z\in \Gamma }{{{b}_{j}}\left( x \right) \delta \left( x-z \right){{\delta }_{ij}}{\rm d}x}} \\
  & =-{{b}_{i}}\left( z \right). \\
\end{aligned}
$$
Therein, the boundary condition vanishes since $p\left( x,0|z,0 \right)=0$ as $\left| {{x}_{j}}-{{z}_{j}} \right|=\varepsilon $. For the second integration, we use integration by parts twice
$$
\begin{aligned}
  & \sum\limits_{k,l=1}^{n}{\int_{x-z\in \Gamma }{\left( {{x}_{i}}-{{z}_{i}} \right)\frac{{{\partial }^{2}}}{\partial {{x}_{k}}\partial {{x}_{l}}}\left[ {{a}_{kl}}\left( x \right)p\left( x,0|z,0 \right) \right]{\rm d}x}} \\
  & =-\sum\limits_{k,l=1}^{n}{\int_{x-z\in \Gamma }{\frac{\partial }{\partial {{x}_{k}}}\left( {{x}_{i}}-{{z}_{i}} \right)\frac{\partial }{\partial {{x}_{l}}} \left[ {{a}_{kl}}\left( x \right)p\left( x,0|z,0 \right) \right]{\rm d}x}} \\
  & =\sum\limits_{k,l=1}^{n}{\int_{x-z\in \Gamma }{\left( {{x}_{i}}-{{z}_{i}} \right){{a}_{kl}}\left( x \right)p\left( x,0|z,0 \right)\frac{\partial {{\delta }_{ik}}}{\partial {{x}_{l}}}{\rm d}x}} \\
  & =0. \\
\end{aligned}
$$
For the third integration, we derive it separately. First, according to Tonelli?s theorem \cite{whitney2012geometric}, we obtain
$$
\begin{aligned}
  & \sum\limits_{j=1}^{n}{\int_{x-z\in \Gamma }{\int_{\mathbbm{R}\backslash \left\{ 0 \right\}}{\left( {{x}_{i}}-{{z}_{i}} \right)p\left( x,0|z,0 \right) W_{j}^{\alpha ,\beta }\left( {y}_{j} \right){\rm d}{y}_{j}}{\rm d}x}} \\
  & =\sum\limits_{j=1}^{n}{\int_{\mathbbm{R}\backslash \left\{ 0 \right\}} {\int_{x-z\in \Gamma }{\left( {{x}_{i}}-{{z}_{i}} \right)\delta \left( x-z \right){\rm d}x}W_{j}^{\alpha ,\beta }\left( {y}_{j} \right){\rm d} {y}_{j}}} \\
  & =0. \\
\end{aligned}
$$
Second,
$$
\begin{aligned}
  & \sum\limits_{j=1}^{n}{\int_{x-z\in \Gamma }{\int_{\mathbbm{R}\backslash \left\{ 0 \right\}}{\left( {{x}_{i}}-{{z}_{i}} \right)p\left( x-{{\sigma }_{j}}(z){y}_{j}{e_{j}},0|z,0 \right)W_{j}^{\alpha ,\beta }\left( {y}_{j} \right){\rm d}{y}_{j}}{\rm d}x}} \\
  & =\int_{x-z\in \Gamma }{\int_{\mathbbm{R}\backslash \left\{ 0 \right\}} {\left( {{x}_{i}}-{{z}_{i}} \right)p\left( x-{{\sigma }_{i}}(z){y}_{i} {e_{i}},0|z,0 \right)W_{i}^{\alpha ,\beta }\left( {y}_{i} \right) {\rm d}{y}_{i}} {\rm d}x} \\
  & =\sigma _{i}^{-1}(z)\int_{-\varepsilon }^{\varepsilon }{{y}_{i}W_{i}^ {\alpha ,\beta }\left( \sigma _{i}^{-1}(z){y}_{i} \right){\rm d}{y}_{i}}. \\
\end{aligned}
$$
Finally, according to Tonelli?s theorem and integration by parts, we have
$$
\begin{aligned}
  & \sum\limits_{j=1}^{n}{\int_{x-z\in \Gamma }{\int_{\mathbbm{R}\backslash \left\{ 0 \right\}}{\left( {{x}_{i}}-{{z}_{i}} \right){{\sigma }_{j}}(z)\chi _{j}^{\alpha }\left( {y}_{j} \right){y}_{j}\frac{\partial }{\partial {{x}_{j}}} p\left( x,0|z,0 \right)W_{j}^{\alpha ,\beta }\left( {y}_{j} \right){\rm d}{y}_{j}}{\rm d}x}} \\
  & =\sum\limits_{j=1}^{n}{\int_{\mathbbm{R}\backslash \left\{ 0 \right\}} {\int_{x-z\in \Gamma }{\left( {{x}_{i}}-{{z}_{i}} \right)\frac{\partial }{\partial {{x}_{j}}}p\left( x,0|z,0 \right){\rm d}x}{{\sigma }_{j}}(z) \chi _{j}^{\alpha }\left( {y}_{j} \right){y}_{j}W_{j}^{\alpha ,\beta }\left( {y}_{j} \right){\rm d}{y}_{j}}} \\
  & =-\sum\limits_{j=1}^{n}{\int_{\mathbbm{R}\backslash \left\{ 0 \right\}} {\int_{x-z\in \Gamma }{p\left( x,0|z,0 \right)\frac{\partial }{\partial {{x}_{j}}}\left( {{x}_{i}}-{{z}_{i}} \right){\rm d}x}{{\sigma }_{j}}(z) \chi _{j}^{\alpha }\left( {y}_{j} \right){y}_{j}W_{j}^{\alpha ,\beta }\left( {y}_{j} \right){\rm d}{y}_{j}}} \\
  & =\left\{ \begin{array}{ll}
   0, & \alpha <1,  \\
   -\sigma _{i}^{-1}(z)\int_{-1}^{1}{{y}_{i}W_{i}^{\alpha ,\beta }\left( \sigma _{i}^{-1}(z){y}_{i} \right){\rm d}{y}_{i}}, & \alpha =1,  \\
   -\sigma _{i}^{-1}(z)\int_{-\infty }^{\infty }{{y}_{i}W_{i}^{\alpha ,\beta }\left( \sigma _{i}^{-1}(z){y}_{i} \right){\rm d}{y}_{i}}, & \alpha >1.  \\
   \end{array} \right. \\
\end{aligned}
$$
Hence,
$$
\underset{t\to 0}{\mathop{\lim }}\,{{t}^{-1}}\int_{x-z\in \Gamma }{\left( {{x}_{i}}-{{z}_{i}} \right)p\left( x,t|z,0 \right){\rm d}x}={{b}_{i}}\left( z \right)+R_{i}^{\alpha ,\beta, \varepsilon}\left( z \right).
$$

\noindent 3) According to the Fokker-Planck equation (\ref{FPK}), we have
$$
\begin{aligned}
  & \underset{t\to 0}{\mathop{\lim }}\,{{t}^{-1}}\int_{x-z\in \Gamma }{\left( {{x}_{i}}-{{z}_{i}} \right)\left( {{x}_{j}}-{{z}_{j}} \right)p\left( x,t|z,0 \right){\rm d}x} \\
  & =\underset{t\to 0}{\mathop{\lim }}\,{{t}^{-1}}\int_{x-z\in \Gamma }{\left( {{x}_{i}}-{{z}_{i}} \right)\left( {{x}_{j}}-{{z}_{j}} \right)\left[ p\left( x,t|z,0 \right)-p\left( x,0|z,0 \right)+p\left( x,0|z,0 \right) \right]{\rm d}x} \\
  & =\int_{x-z\in \Gamma }{\left( {{x}_{i}}-{{z}_{i}} \right)\left( {{x}_{j}}- {{z}_{j}} \right){{\left. \frac{\partial p\left( x,t|z,0 \right)}{\partial t} \right|}_{t=0}}{\rm d}x} \\
  & \ \ \ +\underset{t\to 0}{\mathop{\lim }}\,\left[ {{t}^{-1}}\int_{x-z\in \Gamma }{\left( {{x}_{i}}-{{z}_{i}} \right)\left( {{x}_{j}}-{{z}_{j}} \right)\delta \left( x-z \right){\rm d}x} \right] \\
  & =-\sum\limits_{k=1}^{n}{\int_{x-z\in \Gamma }{\left( {{x}_{i}}-{{z}_{i}} \right)\left( {{x}_{j}}-{{z}_{j}} \right)\frac{\partial }{\partial {{x}_{k}}} \left[ {{b}_{k}}\left( x \right)p\left( x,0|z,0 \right) \right] {\rm d}x}} \\
  & \ \ \ +\frac{1}{2}\sum\limits_{k,l=1}^{n}{\int_{x-z\in \Gamma }{\left( {{x}_{i}}-{{z}_{i}} \right)\left( {{x}_{j}}-{{z}_{j}} \right)\frac{{{\partial }^{2}}}{\partial {{x}_{k}}\partial {{x}_{l}}}\left[ {{a}_{kl}}\left( x \right) p\left( x,0|z,0 \right) \right]{\rm d}x}} \\
  & \ \ \ -\sum\limits_{k=1}^{n}{\int_{x-z\in \Gamma } {\int_{\mathbbm{R} \backslash \left\{ 0 \right\}}{\left( {{x}_{i}}-{{z}_{i}} \right)\left( {{x}_{j}}-{{z}_{j}} \right)\left[ p\left( x,0|z,0 \right)-p\left( x-{{\sigma }_{k}}(z){y}_{k}{e_{k}},0|z,0 \right)-{{\sigma }_{k}}(z)\chi _{k}^{\alpha }\left( {y}_{k} \right){y}_{k}\frac{\partial }{\partial {{x}_{k}}}p\left( x,0|z,0 \right) \right]W_{k}^{\alpha ,\beta }\left( {y}_{k} \right) {\rm d}{y}_{k}} {\rm d}x}}. \\
\end{aligned}
$$
The application of integration by parts into the first term yields
$$
\begin{aligned}
  & \sum\limits_{k=1}^{n}{\int_{x-z\in \Gamma }{\left( {{x}_{i}}-{{z}_{i}} \right) \left( {{x}_{j}}-{{z}_{j}} \right)\frac{\partial }{\partial {{x}_{k}}} \left[ {{b}_{k}}\left( x \right)p\left( x,0|z,0 \right) \right]{\rm d}x}} \\
  & =-\sum\limits_{k=1}^{n}{\int_{x-z\in \Gamma }{{{b}_{k}}\left( x \right) p\left( x,0|z,0 \right)\frac{\partial }{\partial {{x}_{k}}}\left[ \left( {{x}_{i}}-{{z}_{i}} \right)\left( {{x}_{j}}-{{z}_{j}} \right) \right] {\rm d}x}} \\
  & =-\sum\limits_{k=1}^{n}{\int_{x-z\in \Gamma }{{{b}_{k}}\left( x \right) \delta \left( x-z \right)\left[ {{\delta }_{ik}}\left( {{x}_{j}}-{{z}_{j}} \right)+{{\delta }_{jk}}\left( {{x}_{i}}-{{z}_{i}} \right) \right] {\rm d}x}} \\
  & =0. \\
\end{aligned}
$$
Therein, the boundary condition vanishes since $p\left( x,0|z,0 \right)=0$ as $\left| {{x}_{k}}-{{z}_{k}} \right|=\varepsilon $. For the second integration, we use integration by parts again
$$
\begin{aligned}
  & \frac{1}{2}\sum\limits_{k,l=1}^{n}{\int_{x-z\in \Gamma }{\left( {{x}_{i}}- {{z}_{i}} \right)\left( {{x}_{j}}-{{z}_{j}} \right)\frac{{{\partial }^{2}}} {\partial {{x}_{k}}\partial {{x}_{l}}}\left[ {{a}_{kl}}\left( x \right) p\left( x,0|z,0 \right) \right]{\rm d}x}} \\
  & =-\frac{1}{2}\sum\limits_{k,l=1}^{n}{\int_{x-z\in \Gamma }{\frac{\partial }{\partial {{x}_{k}}}\left[ \left( {{x}_{i}}-{{z}_{i}} \right)\left( {{x}_{j}}-{{z}_{j}} \right) \right]\frac{\partial }{\partial {{x}_{l}}}\left[ {{a}_{kl}}\left( x \right)p\left( x,0|z,0 \right) \right]{\rm d}x}} \\
  & =\frac{1}{2}\sum\limits_{k,l=1}^{n}{\int_{x-z\in \Gamma }{{{a}_{kl}}\left( x \right)p\left( x,0|z,0 \right)\frac{\partial }{\partial {{x}_{l}}}\left[ {{\delta }_{ik}}\left( {{x}_{j}}-{{z}_{j}} \right)+{{\delta }_{jk}}\left( {{x}_{i}}-{{z}_{i}} \right) \right]{\rm d}x}} \\
  & =\frac{1}{2}\sum\limits_{k,l=1}^{n}{\int_{x-z\in \Gamma }{{{a}_{kl}}\left( x \right)\delta \left( x-z \right)\left( {{\delta }_{ik}}{{\delta }_{jl}}+{{\delta }_{il}}{{\delta }_{jk}} \right){\rm d}x}} \\
  & =\frac{1}{2}\left[ {{a}_{ij}}\left( z \right)+{{a}_{ji}}\left( z \right) \right] \\
  & ={{a}_{ij}}\left( z \right). \\
\end{aligned}
$$
We still derive the third integration separately. First, according to Tonelli?s theorem, we obtain
$$
\begin{aligned}
  & \sum\limits_{j=1}^{n}{\int_{x-z\in \Gamma }{\int_{\mathbbm{R}\backslash \left\{ 0 \right\}}{\left( {{x}_{i}}-{{z}_{i}} \right)\left( {{x}_{j}}-{{z}_{j}} \right)p\left( x,0|z,0 \right)W_{j}^{\alpha ,\beta }\left( {y}_{j} \right){\rm d}{y}_{j}}{\rm d}x}} \\
  & =\sum\limits_{j=1}^{n}{\int_{\mathbbm{R}\backslash \left\{ 0 \right\}} {\int_{x-z\in \Gamma }{\left( {{x}_{i}}-{{z}_{i}} \right)\left( {{x}_{j}}- {{z}_{j}} \right)\delta \left( x-z \right){\rm d}x}W_{j}^{\alpha ,\beta }\left( {y}_{j} \right){\rm d}{y}_{j}}} \\
  & =0. \\
\end{aligned}
$$
Second, for $i\neq j$,
$$
\sum\limits_{k=1}^{n}{\int_{x-z\in \Gamma }{\int_{\mathbbm{R}\backslash \left\{ 0 \right\}}{\left( {{x}_{i}}-{{z}_{i}} \right)\left( {{x}_{j}}-{{z}_{j}} \right) p\left( x-{{\sigma }_{k}}(z){y}_{k}{e_{k}},0|z,0 \right)W_{k}^{\alpha ,\beta }\left( {y}_{k} \right){\rm d}{y}_{k}}{\rm d}x}}=0
$$
and for $i=j$,
$$
\begin{aligned}
  & \sum\limits_{k=1}^{n}{\int_{x-z\in \Gamma }{\int_{\mathbbm{R}\backslash \left\{ 0 \right\}}{{{\left( {{x}_{i}}-{{z}_{i}} \right)}^{2}}p\left( x-{{\sigma }_{k}}(z){y}_{k}{e_{k}},0|z,0 \right)W_{k}^{\alpha ,\beta }\left( {y}_{k} \right)\textrm{d}{y}_{k}}{\rm d}x}} \\
  & =\int_{x-z\in \Gamma }{\int_{\mathbbm{R}\backslash \left\{ 0 \right\}} {{{\left( {{x}_{i}}-{{z}_{i}} \right)}^{2}}p\left( x-{{\sigma }_{i}}(z) {y}_{i}{e_{i}},0|z,0 \right)W_{i}^{\alpha ,\beta }\left( {y}_{i} \right) \textrm{d}{y}_{i}}{\rm d}x} \\
  & =\sigma _{i}^{-1}(z)\int_{-\varepsilon }^{\varepsilon }{{{y}_{i}^{2}} W_{i}^{\alpha ,\beta }\left( \sigma _{i}^{-1}(z){y}_{i} \right) {\rm d}{y}_{i}}. \\
\end{aligned}
$$
This integration is bounded according to the definition of the jump measure. Finally, according to Tonelli?s theorem and integration by parts, we have
$$
\begin{aligned}
  & \sum\limits_{k=1}^{n}{\int_{x-z\in \Gamma }{\int_{\mathbbm{R}\backslash \left\{ 0 \right\}}{\left( {{x}_{i}}-{{z}_{i}} \right)\left( {{x}_{j}}- {{z}_{j}} \right){{\sigma }_{k}}(z)\chi _{k}^{\alpha }{y}_{k}\frac{\partial }{\partial {{x}_{k}}}p\left( x,0|z,0 \right)W_{k}^{\alpha ,\beta }\left( {y}_{k} \right){\rm d}{y}_{k}}{\rm d}x}} \\
  & =\sum\limits_{k=1}^{n}{\int_{\mathbbm{R}\backslash \left\{ 0 \right\}} {\int_{x-z\in \Gamma }{\left( {{x}_{i}}-{{z}_{i}} \right)\left( {{x}_{j}}- {{z}_{j}} \right)\frac{\partial }{\partial {{x}_{k}}}p\left( x,0|z,0 \right) {\rm d}x}{{\sigma }_{k}}(z)\chi _{k}^{\alpha }{y}_{k}W_{k}^{\alpha ,\beta }\left( {y}_{k}\right){\rm d}{y}_{k}}} \\
  & =-\sum\limits_{k=1}^{n}{\int_{\mathbbm{R}\backslash \left\{ 0 \right\}} {\int_{x-z\in \Gamma }{\delta \left( x-z \right)\left[ {{\delta }_{ik}}\left( {{x}_{j}}-{{z}_{j}} \right)+{{\delta }_{jk}}\left( {{x}_{i}}-{{z}_{i}} \right) \right]{\rm d}x}{{\sigma }_{k}}(z)\chi _{k}^{\alpha }{y}_{k} W_{k}^{\alpha ,\beta }\left( {y}_{k} \right){\rm d}{y}_{k}}} \\
  & =0. \\
\end{aligned}
$$
Hence,
$$
\underset{t\to 0}{\mathop{\lim }}\,{{t}^{-1}}\int_{x-z\in \Gamma }{\left( {{x}_{i}}-{{z}_{i}} \right)\left( {{x}_{j}}-{{z}_{j}} \right)p\left( x,t|z,0 \right){\rm d}x}={{a}_{ij}}\left( z \right)+S_{ij}^{\alpha ,\beta, \varepsilon} \left( z  \right).
$$

\noindent The proof is complete.

\section{Proof of Corollary \ref{cor2}}
\label{App:C}
i)  This is derived directly by integrating the equation in Theorem \ref{thm1}'s first assertion about ${{x}_{i}}-z_i$ on the interval $\left[ {{c}_{1}},\ {{c}_{2}} \right)$.

\noindent ii)  Let the set ${\rm d}U=\left[ {{u}_{1}},\ {{u}_{1}}+{\rm d}{{u}_{1}} \right)\times \left[ {{u}_{2}},\ {{u}_{2}}+{\rm d}{{u}_{2}} \right)\times \cdots \times \left[ {{u}_{n}},\ {{u}_{n}}+{\rm d}{{u}_{n}} \right)$. Then we have
$$
\begin{aligned}
  & \mathbbm{P}\left\{ \left. x\left( t \right)\in {\rm d}U;\ x\left( t \right)-z\in \Gamma  \right|x\left( 0 \right)=z \right\} \\
  & =\mathbbm{P}\left\{ \left. x\left( t \right)\in {\rm d}U \right|x\left( 0 \right)=z;\ x\left( t \right)-z\in \Gamma  \right\}\cdot \mathbbm{P}\left\{ \left. x\left( t \right)-z\in \Gamma  \right|x\left( 0 \right)=z \right\}. \\
\end{aligned}
$$
Thus
$$
\begin{aligned}
  & \int_{x-z\in \Gamma }{\left( {{x}_{i}}-{{z}_{i}} \right)p\left( x,t|z,0 \right){\rm d}x} \\
  & =\int_{u-z\in \Gamma }{\left( {{u}_{i}}-{{z}_{i}} \right)\mathbbm{P}\left\{ \left. x\left( t \right)\in {\rm d}U \right|x\left( 0 \right)=z \right\}} \\
  & =\int_{u-z\in \Gamma }{\left( {{u}_{i}}-{{z}_{i}} \right)\mathbbm{P}\left\{ \left. x\left( t \right)\in {\rm d}U;\ x\left( t \right)-z\in \Gamma  \right|x\left( 0 \right)=z \right\}} \\
  & =\mathbbm{P}\left\{ \left. x\left( t \right)-z\in \Gamma  \right|x\left( 0 \right)=z \right\}\cdot \int_{\left| u-z \right|<\varepsilon }{\left( {{u}_{i}}-{{z}_{i}} \right)\mathbbm{P}\left\{ \left. x\left( t \right)\in {\rm d}U \right|x\left( 0 \right)=z;\ x\left( t \right)-z\in \Gamma  \right\}} \\
  & =\mathbbm{P}\left\{ \left. x\left( t \right)-z\in \Gamma  \right|x\left( 0 \right)=z \right\}\cdot \mathbbm{E}\left[ \left. \left( {{x}_{i}}\left( t \right)-{{z}_{i}} \right) \right|x\left( 0 \right)=z;\ x\left( t \right)-z\in \Gamma  \right]. \\
\end{aligned}
$$
Hence, the conclusion is immediately deduced
$$
\underset{t\to 0}{\mathop{\lim }}\,{{t}^{-1}}\mathbbm{P}\left\{ \left. x\left( t \right)-z\in \Gamma  \right|x\left( 0 \right)=z \right\}\cdot \mathbbm{E}\left[ \left. \left( {{x}_{i}}\left( t \right)-{{z}_{i}} \right) \right|x\left( 0 \right)=z;\ x\left( t \right)-z\in \Gamma  \right]={{b}_{i}}\left( z \right)+R_{i}^{\alpha ,\beta, \varepsilon}\left( z  \right).
$$

\noindent iii) This proof is similar to the second assertion.

\noindent The proof is complete.

\section{Discussion for SDE with rotationally symmetric L\'evy noise}
\label{App:D}
We restate SDE (\ref{sde}) as
\begin{equation} \label{sde2}
	{\rm d}x\left( t \right) =b\left( x\left( t- \right) \right){\rm d}t+ \Lambda \left( x\left( t- \right) \right){\rm d}{B_{t}}+ \sigma(x(t-)) {\rm d}{L_{t}}.
\end{equation}
Here, the L\'evy motion $L_t$ in Eq. (\ref{sde2}) is rotationally symmetric, meaning its jump measure satisfies $\nu\left( {\rm d}\xi \right)= W \left( |\xi| \right){\rm d}\xi$. Additionally, $\sigma(x)$ is a positive scalar function. For this configuration, the corresponding Fokker-Planck equation takes the form
\begin{equation}\label{FPK2}
	\begin{aligned}
		\frac{\partial p}{\partial t}= & -\sum\limits_{i=1}^{n} {\frac{\partial } {\partial {{x}_{i}}} \left[ {{b}_{i}} p\left( x,t|z,0 \right) \right]} +\frac{1}{2}\sum\limits_{i,j=1}^{n} {\frac{{{\partial }^{2}}}{\partial {{x}_{i}}\partial {{x}_{j}}} \left[ {{a}_{ij}}p\left( x,t|z,0 \right) \right]} \\
		& +{\int_{\mathbbm{R}^n\backslash \left\{ 0 \right\}} {\left[ p\left( x+{{\sigma }(x)} y, t|z,0 \right)- p\left( x,t|z,0 \right) \right] W\left( y \right) {\rm d}y}}.
	\end{aligned}
\end{equation}
Therefore, the relationships among the drift, diffusion, L\'evy jump measure, and the transition probability density function $p\left( x,t|z,0 \right)$--solutions to the Fokker-Planck equation (\ref{FPK2})--are established in the following theorem.

\begin{thm} (Relation between stochastic governing law and Fokker-Planck equation)\\
	\label{thm3}
	For every $\varepsilon >0$, the probability density function $p\left( x,t|z,0 \right)$ and the jump measure, drift and diffusion have the following relations:
	\begin{enumerate}[1)]
		\item For every $x$ and $z$ satisfying $\left\| x-z \right\|_2 >\varepsilon $,
		$$
		\underset{t\to 0}{\mathop{\lim }}\,{{t}^{-1}}p \left( x,t|z,0 \right)= \sigma^{-n}(z) W \left( \sigma^{-1}(z) \left( x-z \right) \right).
		$$
		
		\item For $i=1,\ 2,\ \ldots ,\ n$,
		$$
		\underset{t\to 0}{\mathop{\lim }}\,{{t}^{-1}} \int_{\|x-z\|_2 <\varepsilon} {\left( {{x}_{i}}-{{z}_{i}} \right) p\left( x,t|z,0 \right) {\rm d}x} ={{b}_{i}} \left( z \right).
		$$
		
		\item For $i,j=1,\ 2,\ \ldots ,\ n$,
		$$
		\underset{t\to 0}{\mathop{\lim }}\,{{t}^{-1}} \int_{\|x-z\|_2 <\varepsilon} {\left( {{x}_{i}}-{{z}_{i}} \right)\left( {{x}_{j}}-{{z}_{j}} \right) p\left( x,t|z,0 \right){\rm d}x}={{a}_{ij}}\left( z \right)+ S_{ij} \left( z  \right),
		$$
		where $S_{ii} \left( z  \right)=\sigma^{-n}(z) \int_{\|y\|_2 <\varepsilon} {{{y}_{i}^{2}}W\left( \sigma^{-1}(z)y \right){\rm d}y}$ and $S_{ij}\left( z  \right)=0$ for $i\neq j$.
	\end{enumerate}
\end{thm}

This theorem may also be restated as the following nonlocal Kramers-Moyal formulas, which characterize the drift coefficient, diffusion coefficient, and L\'evy jump measure through sample paths--specifically, solutions to the SDE (\ref{sde2}).

\begin{cor}(Nonlocal Kramers-Moyal formulas)\\
	\label{cor4}
	For every $\varepsilon >0$, the sample path solution $x\left( t \right)$ of the stochastic differential equation (\ref{sde2}) and the jump measure, drift and diffusion have the following relations:
	\begin{enumerate}[1)]
		\item For every $m>1$,
		$$
		\underset{t\to 0}{\mathop{\lim }}\,{{t}^{-1}}\mathbbm{P}\left\{ \left. \left\| x \left( t \right)-z \right\|_2 \in \left[ \varepsilon, m\varepsilon \right) \right|x\left( 0 \right)=z \right\} =\sigma ^{-n}(z) \int_{\|y\|_2\in \left[ \varepsilon, m\varepsilon \right)} {W\left( \sigma^{-1}(z) y \right){\rm d}y}.
		$$
		\item For $i=1,\ 2,\ \ldots ,\ n$,
		$$
		\begin{aligned}
			& \underset{t\to 0}{\mathop{\lim }}\,{{t}^{-1}}\mathbbm{P}\left\{ \left. \left\| x \left( t \right)-z \right\|_2 < \varepsilon \right| x\left( 0 \right)=z \right\}\cdot \mathbbm{E}\left[ \left. \left( {{x}_{i}}\left( t \right)-{{z}_{i}} \right) \right| x\left( 0 \right)=z;\ \left\| x \left( t \right)-z \right\|_2 < \varepsilon  \right] \\
			& ={{b}_{i}}\left( z \right). \\
		\end{aligned}
		$$
		\item For $i,j=1,\ 2,\ \ldots ,\ n$,
		$$
		\begin{aligned}
			& \underset{t\to 0}{\mathop{\lim }}\,{{t}^{-1}}\mathbbm{P}\left\{ \left. \left\| x \left( t \right)-z \right\|_2 < \varepsilon  \right|x\left( 0 \right)=z \right\}\cdot \mathbbm{E}\left[ \left. \left( {{x}_{i}}\left( t \right)-{{z}_{i}} \right)\left( {{x}_{j}}\left( t \right)-{{z}_{j}} \right) \right|x\left( 0 \right)=z;\ \left\| x \left( t \right)-z \right\|_2 < \varepsilon  \right] \\
			& ={{a}_{ij}}\left( z \right)+S_{ij} \left( z  \right). \\
		\end{aligned}
		$$
	\end{enumerate}
\end{cor}

The proofs of Theorem \ref{thm3} and Corollary \ref{cor4} follow straightforwardly by combining the arguments in Ref. \cite{YangLi2020a} with the proofs of Theorem \ref{thm1} and Corollary \ref{cor2}; they are therefore omitted.

The data-driven algorithms for discovering the drift and diffusion terms of SDE (\ref{sde2}) are similar to those in Section \ref{DDAsec}, with the key modification being the change of the integral domain from a cube to a ball. The primary distinction lies in the process of learning the L\'evy jump measure and noise intensity. To avoid the curse of dimensionality, we assume that the noise intensity function also possesses rotational symmetry, i.e., $\sigma(x)=\sigma(r)$ where $r=\left\| x \right\|_2$.

Similarly, We construct a new dataset ${\mathbbm{Y}}= \left[ y^{(1)},\ y^{(2)},\ \cdots ,\ y^{(M)} \right]$ where ${{y}^{(j)}}= \|{{x}^{(j)}}- {{z}^{(j)}}\|_2$ for $j=1,\ 2,\ \ldots ,\ M$. Since Corollary \ref{cor4}'s first assertion depends on $\|z\|_2$, we partition its phase space $0 \leq \|z\|_2 <z^{\text{max}}$ into $N_c$ intervals: $\left[ 0,\ \Delta z \right)$, $\left[ \Delta z,\ 2\Delta z \right)$, $\ldots$, $\left[ z^{\text{max}}- \Delta z,\ z^{\text{max}} \right)$, where $\Delta z= z^{\text{max}}/N_c$. Using the midpoint $\hat{z}_{l}= (l-0.5) \Delta z$, $l=1,\ 2,\ \ldots ,\ N_c$ of each interval, we approximate ${{\sigma}}(r)$ as the constant ${{\sigma }}(\hat{z}_{l})$. For each interval $\|z\|_2 \in \left[ (l-1) \Delta z,\ l\Delta z \right)$, let $\mathbbm{Z}_{l}$ denote the corresponding subset of $\mathbbm{Z}$, with $\mathbbm{X}_{l}$ and $\mathbbm{Y}_{l}$ being the associated $M_l$-element datasets. We then estimate the conditional probability as the proportion of points in ${{Y}_{l}}$ falling within $\left[ \varepsilon ,\ m\varepsilon  \right)$ relative to $M_l$. We choose $N+1$ intervals $\left[ \varepsilon ,\ m\varepsilon  \right)$, $\left[ m\varepsilon ,\ {{m}^{2}}\varepsilon  \right)$, $\ldots$, $\left[ {{m}^{N}}\varepsilon ,\ {{m}^{N+1}}\varepsilon  \right)$ using the positive integer $N$, the positive real number $\varepsilon $, and the real number $m>1$. Let $n_{0,l},\ n_{1,l},\ \ldots ,\ n_{N,l}$ represent the number of points from dataset ${{Y}_{l}}$ falling into these intervals. 

From Corollary \ref{cor4}'s first assertion:
\begin{equation}\label{Aeq2}
\begin{aligned}
	& {{h}^{-1}}\mathbbm{P}\left\{ \left. \| {{x}}\left( h \right)-z \|_2 \in \left[ {{m}^{k}}\varepsilon ,\ {{m}^{k+1}}\varepsilon  \right) \right|x\left( 0 \right)=z \right\}\approx {{h}^{-1}}{{M}_l^{-1}}n_{k,l}, \\
	& k=0,\ 1,\ \ldots ,\ N, \ l=1,\ 2,\ \ldots,\ N_c.
\end{aligned}
\end{equation}
Consider $n=2$ for concreteness. Then its right-hand side integration can be computed as
\begin{equation}\label{Aeq3}
\begin{aligned}
	& {{\sigma }^{-n} (\hat{z}_{l})} \int_{\|y\|_2 \in \left[ {{m}^{k}}\varepsilon ,\ {{m}^{k+1}}\varepsilon  \right) } {W\left( {{\sigma }^{-1} (\hat{z}_{l})} y \right) {\rm d}y} \\
	& ={{\sigma }^{\alpha} (\hat{z}_{l})} c(2,\alpha) \int_{\|y\|_2 \in \left[ {{m}^{k}}\varepsilon ,\ {{m}^{k+1}}\varepsilon  \right) } {{\|y\|_2 ^{-\left( 2+\alpha  \right)}}{\rm d}y} \\
	& ={{\sigma }^{\alpha} (\hat{z}_{l})} c(2,\alpha) \int_0^{2\pi} {\rm d} \theta \int_{{m}^{k}\varepsilon}^{{{m}^{k+1}}\varepsilon} {{r ^{-\left( 1+\alpha  \right)}}{\rm d}r} \\
	& =2\pi {{\sigma }^{\alpha} (\hat{z}_{l})} c(2,\alpha) {{\alpha }^{-1}}{{\varepsilon }^{-\alpha }}{{m}^{-k\alpha }}\left( 1-{{m}^{-\alpha }} \right), \\
	& k=0,\ 1,\ \ldots ,\ N, \ l=1,\ 2,\ \ldots,\ N_c,
\end{aligned}
\end{equation}
where $c(n,\alpha)= \frac{\alpha \Gamma((n+\alpha)/2)} {2^{1-\alpha} \pi^{n/2} \Gamma(1-\alpha/2)}$. Combining these two equations yields $(N+1)N_c$ equalities
\begin{equation}\label{Aeq4}
	\begin{aligned}
		& 2\pi {{\sigma }^{\alpha} (\hat{z}_{l})} c(2,\alpha) {{\alpha }^{-1}}{{\varepsilon }^{-\alpha }}{{m}^{-k\alpha }}\left( 1-{{m}^{-\alpha }} \right)={{h}^{-1}}{{M}_l^{-1}}n_{k,l}, \\
		& k=0,\ 1,\ \ldots ,\ N, \ l=1,\ 2,\ \ldots,\ N_c.
	\end{aligned}
\end{equation}
Summing Eq. (\ref{Aeq4}) over index $l$ gives
\begin{equation}\label{Aeq5}
	\begin{aligned}
		& 2\pi c(2,\alpha) {{{\alpha }^{-1}}{{\varepsilon }^{-\alpha }} {{m}^{-k\alpha }}\left( 1-{{m}^{-\alpha }} \right)} \sum_{l=1}^{N_c} {{\sigma }^{\alpha} (\hat{z}_{l})} ={{h}^{-1}}\sum_{l=1}^{N_c} {{M}_l^{-1}} n_{k,l}, \\
		& k=0,\ 1,\ \ldots ,\ N.
	\end{aligned}
\end{equation}
Taking ratios of the first equation ($k=0$) to the remaining $N$ equations yields
\begin{equation}\label{Aeq6}
	\alpha ={{\left( k\ln m \right)}^{-1}}\ln \frac{\sum_{l=1}^{N_c} {{M}_l^{-1}} n_{0,l}}{\sum_{l=1}^{N_c} {{M}_l^{-1}} n_{k,l}},\quad k=1,\ 2,\ \ldots ,\ N.
\end{equation}
For $N=1$, this provides the optimal estimate $\tilde{\alpha}$. When $N \geq 2$, we take $\tilde{\alpha}$ as the mean value of Eqs. (\ref{Aeq6}).

Summing Eq. (\ref{Aeq4}) over index $k$ yields
\begin{equation}\label{Aeq7}
	\begin{aligned}
		& {2\pi {{\sigma }^{\alpha} (\hat{z}_{l})} c(2,\alpha) {{\alpha }^{-1}} {{\varepsilon }^{-\alpha }}\left( 1-{{m}^{-(N+1)\alpha }} \right)}= {{h}^{-1}} {{M}_l^{-1}} \sum_{k=0}^{N}n_{k,l}, \\
		& l=1,\ 2,\ \ldots,\ N_c.
	\end{aligned}
\end{equation}
Thus, the noise intensity ${\sigma } (\hat{z}_{l})$ can be approximated as
\begin{equation}\label{Aeq8}
	{\sigma } (\hat{z}_{l}) ={{\left[ \frac{\tilde{\alpha }{{\varepsilon }^{{\tilde{\alpha }}}} \sum_{k=0}^{N}n_{k,l}} {2\pi c(2,\alpha) hM_l \left( 1-{{m}^{-(N+1)\alpha }} \right)} \right]} ^{{1}/{{\tilde{\alpha }}}\;}},\quad l=1,\ 2,\ \ldots,\ N_c.
\end{equation}
For multi-dimensional systems, the sole modification occurs in Eq. (\ref{Aeq8}) concerning the constant $2\pi c(2,\alpha)$, as specified by Eq. (\ref{Aeq3}). Define $\mathbbm{Q}= \{ {\tilde{\sigma }}(\hat{z}_{l}): l=1, 2, \ldots,\ N_c \}$. Regression techniques such as sparse regression or neural networks can also be employed to approximate the noise intensity function $\sigma_{\theta}(\|x\|_2)$ by minimizing the loss function
\begin{equation}\label{Aloss}
	\mathcal{L}= \frac{1}{N_c} \sum_{l=1}^{N_c} | \sigma_{\theta}(\hat{z}_{l})- {\tilde{\sigma }}(\hat{z}_{l}) |^2.
\end{equation}

\section{Sparse identification of nonlinear dynamics (SINDy)}
\label{App:E}
The Sparse Identification of Nonlinear Dynamics (SINDy) method, proposed by Brunton et al. \cite{SINDy1}, addresses the sparse regression problem for deterministic systems
\begin{equation}\label{sindy}
	\tilde{c} = \argmin_c (\left\| B-Ac \right\|_2^2 +\lambda \left\| c \right\|_1).
\end{equation}
The algorithm proceeds as follows:
\begin{itemize}
	\item Solve the unconstrained minimization $\hat{c} = \argmin_c \left\| B-Ac \right\|_2^2$ via the normal equation $\hat{c} = (A^TA)^{-1} A^TB$, yielding a non-sparse solution $\hat{c}$.
	
	\item Set coefficients in $\hat{c}$ smaller than a predefined threshold $\lambda$ to zero.
	
	\item Iterate steps 1?2 on the remaining coefficients until no values fall below $\lambda$.
\end{itemize}
The parameter $\lambda$ controls solution sparsity and requires careful calibration.

\section{Binning}
\label{App:F}
For systems of relatively low dimension (e.g., $n=1, 2, 3$), binning operations serve as data pre-processing for learning drift and diffusion coefficients, reducing error prior to stepwise sparse regression. The core methodology involves:
\begin{itemize}
	\item Discretizing phase space into $Q$ bins.
	
	\item Transforming dataset $\mathbbm{\hat{Z}}$ into $\mathbbm{\bar{Z}} = \{\bar{z}^{(j)}: j=1,2, ..., Q\}$ and $\mathbbm{W} = \{w^{(j)}: j=1,2, ..., Q\}$, where $\bar{z}^{(j)}$ is the $j$-bin centroid and $w^{(j)}$ represents the bin's data weight ratio.
	
	\item Constructing matrix $A_Q$ via candidate dictionary and $\mathbbm{\bar{Z}}$.
	
	\item Defining weight matrix $W_Q= \text{diag} \{ \mathbbm{W} \}$.
	
	\item Computing $B_Q$ by averaging matrix $B$ entries within bins.
\end{itemize}
This transforms the sparse regression in Eq. (\ref{sindy}) to
\begin{equation}\label{binning}
	\tilde{c} = \argmin_c (\left\| W_Q B_Q- W_Q A_Q c \right\|_2^2 +\lambda \left\| c \right\|_1).
\end{equation}

\section{Stepwise sparse regressor and cross-validation}
\label{App:G}
Boninsegna et al. \cite{Boninsegna2018} proposed the stepwise sparse regressor and cross-validation algorithms for solving the minimization problem (\ref{sindy}), which constitute a modification of the SINDy method. The stepwise sparse regressor is an iterative algorithm comprising:
\begin{itemize}
	\item Solving the unconstrained least squares problem $\tilde{c} = \argmin_c \left\| B-Ac \right\|_2^2$ via the normal equation to obtain a non-sparse solution $\tilde{c}$.
	
	\item Eliminating the smallest-magnitude coefficient in $\tilde{c}$ and its corresponding function
	\begin{equation}\label{ssrstep2}
		\tilde{c}_i=0: i=\argmin_k \| \tilde{c}_k \|.
	\end{equation}
	
	\item Recomputing the unconstrained regression with remaining functions.
	
	\item Iterating until achieving a satisfactory cross-validation score.
\end{itemize}

The sparse solution from the stepwise sparse regressor algorithm is expressed as $\tilde{c} = SSR_q (A,B)$, where $q$ specifies the sparsity degree (i.e., the number of zero coefficients). Cross-validation is then performed on this $q$-sparse solution $SSR_q (A,B)$, with its score denoted $\delta[SSR_q]$.

Cross-validation partitions the dataset into multiple segments, successively using each segment as the test set to evaluate models trained on the remaining data. This balances algorithmic accuracy and model simplicity. Specifically:
\begin{itemize}
	\item Randomly divide the full dataset $\mathbbm{D}$ into $K_{cv}$ subsets: $\cup_i \mathbbm{D}_i =\mathbbm{D}$, $\mathbbm{D}_i \cap \mathbbm{D}_j =\emptyset$.
	
	\item Let $A_i$ and $B_i$ denote submatrices of $A$ and $B$ corresponding to $\mathbbm{D}_i$ where $i=1,2,...,K_{cv}$.
	
	\item Let $A_{-i}$ and $B_{-i}$ denote submatrices with $\mathbbm{D}_i$-associated rows omitted where $i=1,2,...,K_{cv}$. The cross-validation score for the sparse solution is then defined as:
	\begin{equation}\label{cvscore}
		\delta^2[SSR_q]= \frac{1}{K_{cv}} \sum\limits_{i=1}^{K_{cv}} \| B_i- A_i \cdot SSR_q (A_{-i},B_{-i}) \|_2^2.
	\end{equation}
\end{itemize}
Thus $\delta[SSR_q]$ is a function of $q$. The optimal sparse solution occurs at an abrupt transition point on the score curve, satisfying
\begin{equation}\label{cvcond}
	\frac{\delta_[\tilde{q}]}{\delta_[\tilde{q}-1]} \approx 1, \ \frac{\delta_[\tilde{q}+1]}{\delta_[\tilde{q}]} \gg 1.
\end{equation}

\section{Proof of Theorem \ref{thmerror1}}
\label{App:H}
\noindent 1) According to Corollary \ref{cor2}' first assertion, the short time transition probability is approximated as
\begin{equation*}
	\mathbbm{P}\left\{ \left. x \left( t \right)-z \in \left[ {{c}_{1}}, {{c}_{2}} \right) \right|x\left( 0 \right)=z \right\}= h \cdot \sigma^{-1}(z) \int_{{{c}_{1}}}^{{{c}_{2}}} {W^{\alpha ,\beta }\left( \sigma^{-1}(z) y \right){\rm d}y} + O(h^2).
\end{equation*}
For the intervals $\left[ {{m}^{k}}\varepsilon,\ {{m}^{k+1}}\varepsilon  \right)$ and $\left[ -{{m}^{k+1}}\varepsilon,\  -{{m}^{k}}\varepsilon  \right)$, define
\begin{equation*}
	\begin{aligned}
	& p_k^+(z)= \sigma^{-1}(z) \int_{{{m}^{k}}\varepsilon} ^{{{m}^{k+1}}\varepsilon} {W^{\alpha ,\beta }\left( \sigma^{-1}(z) y \right){\rm d}y}, \\
	& p_k^-(z)= \sigma^{-1}(z) \int_{-{{m}^{k+1}}\varepsilon} ^{-{{m}^{k}}\varepsilon} {W^{\alpha ,\beta }\left( \sigma^{-1}(z) y \right){\rm d}y}.
	\end{aligned}
\end{equation*}

Within the subinterval corresponding to $\hat{z}_l$, count statistics $n_{k,l}^{\pm}$ obey binomial distribution:
\begin{equation*}
	n_{k,l}^{\pm} \sim \text{Binomial} (M_l, h p_k^{\pm}(\hat{z}_l) + O(h^2) + O(hN_c^{-1})).
\end{equation*}
According to law of large numbers, we have
\begin{equation*}
	\frac{n_{k,l}^{\pm}}{M_l} = h p_k^{\pm}(\hat{z}_l) + O(h^2) + O(hN_c^{-1}) + O_{\mathbbm{P}} \left( \sqrt{\frac{h p_k^{\pm}(\hat{z}_l)}{M_l}} \right).
\end{equation*}

Denote $p_k= p_k^{+}(\hat{z}_l) + p_k^{-}(\hat{z}_l)$. Then we have
\begin{equation*}
	\begin{aligned}
	\tilde{\alpha}_k &= {{\left( k\ln m \right)}^{-1}}\ln \frac{\sum_{l=1}^{N_c} {{M}_l^{-1}} (n_{0,l}^{+}+n_{0,l}^{-})}{\sum_{l=1}^{N_c} {{M}_l^{-1}} (n_{k,l}^{+}+n_{k,l}^{-})} \\
	&= {{\left( k\ln m \right)}^{-1}}\ln \frac{\sum_{l=1}^{N_c} \left[ h p_0^{+}(\hat{z}_l) +h p_0^{-}(\hat{z}_l) + O(h^2) + O(hN_c^{-1}) + O_{\mathbbm{P}} \left( \sqrt{\frac{h p_0^{+}(\hat{z}_l)}{M_l}} \right) + O_{\mathbbm{P}} \left( \sqrt{\frac{h p_0^{-}(\hat{z}_l)}{M_l}} \right) \right] } {\sum_{l=1}^{N_c} \left[ h p_k^{+}(\hat{z}_l) +h p_k^{-}(\hat{z}_l) + O(h^2) + O(hN_c^{-1}) + O_{\mathbbm{P}} \left( \sqrt{\frac{h p_k^{+}(\hat{z}_l)}{M_l}} \right) + O_{\mathbbm{P}} \left( \sqrt{\frac{h p_k^{-}(\hat{z}_l)}{M_l}} \right) \right] } \\
	&= {{\left( k\ln m \right)}^{-1}}\ln \frac{\sum_{l=1}^{N_c} \left[ p_0 + O(h) + O(N_c^{-1}) + O_{\mathbbm{P}} \left( \sqrt{M^{-1}h^{-1}N_c} \right) \right] } {\sum_{l=1}^{N_c} p_k } \\
	&= {{\left( k\ln m \right)}^{-1}}\ln \frac{ p_0 + O(h) + O(N_c^{-1}) + O_{\mathbbm{P}} \left( \sqrt{M^{-1}h^{-1}N_c} \right) } { p_k } \\
	\end{aligned}
\end{equation*}
Due to the fact
\begin{equation*}
	\begin{aligned}
		p_k^+(z) &= \sigma^{-1}(z) \int_{{{m}^{k}}\varepsilon} ^{{{m}^{k+1}}\varepsilon} {W^{\alpha ,\beta }\left( \sigma^{-1}(z) y \right){\rm d}y} \\
		&= {{{\sigma }^{\alpha }} (z) {{k}_{\alpha }}\left( 1+\beta  \right){{\alpha }^{-1}}{{\varepsilon }^{-\alpha }}{{m}^{-k\alpha }}\left( 1-{{m}^{-\alpha }} \right)}/{2},  \\
		p_k^-(z) &= \sigma^{-1}(z) \int_{-{{m}^{k+1}}\varepsilon} ^{-{{m}^{k}}\varepsilon} {W^{\alpha ,\beta }\left( \sigma^{-1}(z) y \right){\rm d}y}  \\
		&= {{{\sigma }^{\alpha }} (z) {{k}_{\alpha }}\left( 1-\beta  \right){{\alpha }^{-1}}{{\varepsilon }^{-\alpha }}{{m}^{-k\alpha }}\left( 1-{{m}^{-\alpha }} \right)}/{2},
	\end{aligned}
\end{equation*}
the true parameter $\alpha$ satisfies
\begin{equation*}
\alpha= {{\left( k\ln m \right)}^{-1}}\ln \frac{ p_0 } { p_k }.
\end{equation*}
Therefore, we obtain
\begin{equation*}
	\begin{aligned}
	|\tilde{\alpha}_k-\alpha| &= {{\left( k\ln m \right)}^{-1}}\ln \frac{ p_0 + O(h) + O(N_c^{-1}) + O_{\mathbbm{P}} \left( \sqrt{M^{-1}h^{-1}N_c} \right) } { p_0 } \\
	&= O(h) + O(N_c^{-1}) + O_{\mathbbm{P}} \left( \sqrt{M^{-1}h^{-1}N_c} \right),
	\end{aligned}
\end{equation*}
and the rate will remained unchanged after averaging about $k$.

\noindent 2) According to Eq. (\ref{eq7a}), we have
\begin{equation*}
	\begin{aligned}
	\rho &= \frac{\sum_{l=1}^{N_c} \sum_{k=0}^{N} {{M}_l^{-1}} {n_{k}^{-}}} {\sum_{l=1}^{N_c} \sum_{k=0}^{N} {{M}_l^{-1}} {n_{k}^{+}}} \\
	&= \frac{\sum_{l=1}^{N_c} \sum_{k=0}^{N} \left[ h p_k^{-}(\hat{z}_l) + O(h^2) + O(hN_c^{-1}) + O_{\mathbbm{P}} \left( \sqrt{\frac{h p_k^{-}(\hat{z}_l)}{M_l}} \right) \right] } {\sum_{l=1}^{N_c} \sum_{k=0}^{N} \left[ h p_k^{+}(\hat{z}_l) + O(h^2) + O(hN_c^{-1}) + O_{\mathbbm{P}} \left( \sqrt{\frac{h p_k^{+}(\hat{z}_l)}{M_l}} \right) \right]} \\
	&= \frac{\sum_{l=1}^{N_c} \sum_{k=0}^{N} \left[ p_k^{-} + O(h) + O(N_c^{-1}) + O_{\mathbbm{P}} \left( \sqrt{\frac{p_k^{-}}{M_l h}} \right) \right] } {\sum_{l=1}^{N_c} \sum_{k=0}^{N} \left[ p_k^{+} + O(h) + O(N_c^{-1}) + O_{\mathbbm{P}} \left( \sqrt{\frac{p_k^{+}}{M_l h}} \right) \right]}.
    \end{aligned}
\end{equation*}
Then
\begin{equation*}
	\begin{aligned}
		\tilde{\beta} &= \frac{1-\rho} {1+\rho} \\
		&= \frac{\sum_{l=1}^{N_c} \sum_{k=0}^{N} \left[ p_k^{+} + O(h) + O(N_c^{-1}) + O_{\mathbbm{P}} \left( \sqrt{\frac{p_k^{+}}{M_l h}} \right) - p_k^{-} - O(h) - O(N_c^{-1}) - O_{\mathbbm{P}} \left( \sqrt{\frac{p_k^{-}}{M_l h}} \right) \right] } {\sum_{l=1}^{N_c} \sum_{k=0}^{N} p_k} \\
		&= \frac{\sum_{l=1}^{N_c} \sum_{k=0}^{N} \left[ p_k^{+} - p_k^{-} + O(h) + O(N_c^{-1}) + O_{\mathbbm{P}} \left( \sqrt{M^{-1}h^{-1}N_c} \right) \right] } {\sum_{l=1}^{N_c} \sum_{k=0}^{N} p_k}.
	\end{aligned}
\end{equation*}
Thus
\begin{equation*}
	\begin{aligned}
		|\tilde{\beta}-\beta| &= \frac{\sum_{l=1}^{N_c} \sum_{k=0}^{N} \left[ O(h) + O(N_c^{-1}) + O_{\mathbbm{P}} \left( \sqrt{M^{-1}h^{-1}N_c} \right) \right] } {\sum_{l=1}^{N_c} \sum_{k=0}^{N} p_k} \\
		&= O(h) + O(N_c^{-1}) + O_{\mathbbm{P}} \left( \sqrt{M^{-1}h^{-1}N_c} \right).
	\end{aligned}
\end{equation*}

\noindent 3) From Eq. (\ref{eq8}), we derive
\begin{equation*}
	\begin{aligned}
		\tilde{\sigma}_l^{\alpha} &= \frac{\tilde{\alpha }{{\varepsilon }^{{\tilde{\alpha }}}} \sum_{k=0}^{N}(n_{k,l}^{+}+n_{k,l}^{-})} {{{k}_{{ \tilde{\alpha }}}} hM_l \left( 1-{{m}^{-(N+1)\tilde{\alpha } }} \right)} \\
		&= \frac{\alpha{{\varepsilon }^{\alpha}} \sum_{k=0}^{N}(n_{k,l}^{+}+n_{k,l}^{-})} {{{k}_{\alpha}} hM_l \left( 1-{{m}^{-(N+1)\alpha }} \right)} \left[ 1+ O_{\mathbbm{P}} (|\tilde{\alpha}-\alpha|) \right] \\
		&= \frac{\alpha{{\varepsilon }^{\alpha}} \sum_{k=0}^{N} \left[ p_k + O(h) + O(N_c^{-1}) + O_{\mathbbm{P}} \left( \sqrt{M^{-1}h^{-1}N_c} \right) \right]} {{{k}_{\alpha}} \left( 1-{{m}^{-(N+1)\alpha }} \right)} \left[ 1 + O(h) + O(N_c^{-1}) + O_{\mathbbm{P}} \left( \sqrt{M^{-1}h^{-1}N_c} \right) \right] \\
		&= \frac{\alpha{{\varepsilon }^{\alpha}} \sum_{k=0}^{N} \left[ p_k + O(h) + O(N_c^{-1}) + O_{\mathbbm{P}} \left( \sqrt{M^{-1}h^{-1}N_c} \right) \right]} {{{k}_{\alpha}} \left( 1-{{m}^{-(N+1)\alpha }} \right)}.
	\end{aligned}
\end{equation*}
Thus
\begin{equation*}
	|\tilde{\sigma}_l^{\alpha} - \sigma_l^{\alpha}| = O(h) + O(N_c^{-1}) + O_{\mathbbm{P}} \left( \sqrt{M^{-1}h^{-1}N_c} \right).
\end{equation*}
According to Taylor expansion, we obtain
\begin{equation*}
	|\tilde{\sigma}(\hat{z}_l) - \sigma(\hat{z}_l)| = O(h) + O(N_c^{-1}) + O_{\mathbbm{P}} \left( \sqrt{M^{-1}h^{-1}N_c} \right),
\end{equation*}
and the rate will remained unchanged after averaging about $l$.

\section{Proof of Theorem \ref{thmerror2}}
\label{App:I}
In this one-dimensional case, the set $\Gamma= [-\varepsilon, \varepsilon]$. Define
\begin{equation*}
    p(z)= \sigma^{-1}(z) \int_{y \notin \Gamma} {W^{\alpha ,\beta }\left( \sigma^{-1}(z) y \right){\rm d}y}.
\end{equation*}
Then the conditional probabilities have the following approximations
\begin{equation*}
	\begin{aligned}
	& \mathbbm{P}\left\{ \left. x \left( h \right)-z \notin \Gamma \right|x\left( 0 \right)=z \right\}= h p(z) + o(h), \\
	& \mathbbm{P}\left\{ \left. x \left( h \right)-z \in \Gamma \right|x\left( 0 \right)=z \right\}= 1- h p(z) + o(h).
	\end{aligned}
\end{equation*}
Denote $\mathcal{D}$ as the domain interested in phase space and $\rho_z$ as the measure of the initial distribution of $z$. Then we have
\begin{equation*}
	\begin{aligned}
		& \int_{\mathcal{D}} \mathbbm{P}\left\{ \left. x \left( h \right)-z \in \Gamma \right|x\left( 0 \right)=z \right\} \rho_z({\rm d}z) \\
		& = 1- h \int_{\mathcal{D}} p(z) \rho_z({\rm d}z) + o(h) \\
		& \triangleq 1- h \bar{p} + o(h).
	\end{aligned}
\end{equation*}
On the other hand, we also obtain
\begin{equation*}
		\int_{\mathcal{D}} \mathbbm{P}\left\{ \left. x \left( h \right)-z \in \Gamma \right|x\left( 0 \right)=z \right\} \rho_z({\rm d}z) = \hat{M}M^{-1} + O_{\mathbbm{P}} (\sqrt{M^{-1}h}).
\end{equation*}
Therefore, we derive
\begin{equation*}
	\hat{M}M^{-1} = \mathbbm{P}\left\{ \left. x \left( h \right)-z \in \Gamma \right|x\left( 0 \right)=z \right\} + h(p(z)-\bar{p}) + O_{\mathbbm{P}} (\sqrt{M^{-1}h}) +o(h).
\end{equation*}

Note that the regression problems (\ref{converg1}) and (\ref{converg2}) are equivalent to the normal equations
\begin{equation*}
	\begin{aligned}
		& A^TAc=A^TB_b, \\
		& A^TAd=A^TB_a.
	\end{aligned}
\end{equation*}
Given ergodicity of the random process, we have
\begin{equation*}
	\begin{aligned}
		\frac{1}{\hat{M}} [A^TA]_{i,j} & = \frac{1}{\hat{M}} \sum_{k=1}^{\hat{M}} \psi_i (\hat{z}^{(k)}) \psi_j (\hat{z}^{(k)}) \\
		& \xrightarrow{M \to \infty} \int_{\mathcal{D}} \psi_i (z) \psi_j (z) \mu({\rm d}z), \\
		\frac{1}{\hat{M}} [A^T B_b]_i & = \frac{1}{\hat{M}} \sum_{k=1}^{\hat{M}} \psi_i (\hat{z}^{(k)}) [\hat{M}M^{-1}h^{-1} (\hat{x}^{(k)} - \hat{z}^{(k)}) - R^{\alpha, \beta, \varepsilon} (\hat{z}^{(k)})] \\
		& \xrightarrow{M \to \infty} \int_{\mathcal{D}} \psi_i (z) \{ [\mathbbm{P}\left\{ \left. x \left( h \right)-z \in \Gamma \right|x\left( 0 \right)=z \right\} + O(h)] \cdot \\
		& \ \ \ \ \ \ \ \ \ \ \ \ h^{-1} \mathbbm{E}^z [x(h)-z | x(h)-z \in \Gamma] - R^{\alpha, \beta, \varepsilon} (z) \} \mu({\rm d}z), \\
		\frac{1}{\hat{M}} [A^T B_a]_i & = \frac{1}{\hat{M}} \sum_{k=1}^{\hat{M}} \psi_i (\hat{z}^{(k)}) [\hat{M}M^{-1}h^{-1} (\hat{x}^{(k)} - \hat{z}^{(k)})^2 - S^{\alpha, \beta, \varepsilon} (\hat{z}^{(k)})] \\
		& \xrightarrow{M \to \infty} \int_{\mathcal{D}} \psi_i (z) \{ [\mathbbm{P}\left\{ \left. x \left( h \right)-z \in \Gamma \right|x\left( 0 \right)=z \right\} + O(h)] \cdot \\
		&\ \ \ \ \ \ \ \ \ \ \ \  h^{-1} \mathbbm{E}^z [(x(h)-z)^2 | x(h)-z \in \Gamma] - S^{\alpha, \beta, \varepsilon} (z) \} \mu({\rm d}z).
	\end{aligned}
\end{equation*}
In the limit $h \to 0$, the nonlocal Kramers-Moyal formulas imply that
\begin{equation*}
	\begin{aligned}
		& \int_{\mathcal{D}} \psi_i (z) \{ [\mathbbm{P}\left\{ \left. x \left( h \right)-z \in \Gamma \right|x\left( 0 \right)=z \right\} + O(h)] \cdot  \\
		& h^{-1} \mathbbm{E}^z [x(h)-z | x(h)-z \in \Gamma] - R^{\alpha, \beta, \varepsilon} (z) \} \mu({\rm d}z) \\
		& \xrightarrow{h \to 0} \int_{\mathcal{D}} \psi_i (z) b(z) \mu({\rm d}z), \\
		& \int_{\mathcal{D}} \psi_i (z) \{ [\mathbbm{P}\left\{ \left. x \left( h \right)-z \in \Gamma \right|x\left( 0 \right)=z \right\} + O(h)] \cdot \\
		& h^{-1} \mathbbm{E}^z [(x(h)-z)^2 | x(h)-z \in \Gamma] - S^{\alpha, \beta, \varepsilon} (z) \} \mu({\rm d}z) \\
		& \xrightarrow{h \to 0} \int_{\mathcal{D}} \psi_i (z) a(z) \mu({\rm d}z).
	\end{aligned}
\end{equation*}
These confirm that the solutions of the regression problems (\ref{converg1}) and (\ref{converg2}) converge to the solutions of the best approximations for drift $b$ and diffusion $a$ in the space $L_{\mu}^2$.


%
%


\bibliographystyle{abbrv}
\bibliography{manuscript}

\end{document}